\documentclass[a4paper,11pt,final]{article}

\usepackage{ulem} 
\usepackage{graphicx,subcaption,booktabs,paralist,multirow,framed}
\usepackage[bottom]{footmisc} 
\usepackage[authoryear,round]{natbib}
\usepackage[usenames,dvipsnames]{xcolor}
\usepackage[colorlinks=true,citecolor=Mahogany,linkcolor=Mahogany,urlcolor=Mahogany]{hyperref}
\usepackage{doi}
\usepackage{mathtools}
\usepackage{bm}
\usepackage{ushort}
\usepackage{xspace} 
\usepackage{mathpazo} 
\usepackage[margin=1.in]{geometry}
\usepackage{setspace}
\usepackage{pifont} 
\usepackage{nicefrac} 
\usepackage{stmaryrd}
\usepackage{amsthm}
\usepackage{float}
\definecolor{c1}{HTML}{4477AA}
\definecolor{c2}{HTML}{CC6677}

\usepackage[obeyDraft]{todonotes}
\usepackage{datetime}
\newdateformat{monthyeardate}{\monthname[\THEMONTH], \THEYEAR}
\usepackage{tabularx}

\newcolumntype{L}{>{\centering\arraybackslash}X}
\usepackage{siunitx}
\begingroup
\catcode`\_=11
\catcode`\:=11
 
\endgroup

\newcommand{\Rm}{{\mathbb R}}
\newcommand{\Em}{{\mathbb E}}
\newcommand{\be}{\begin{equation}}
\newcommand{\ee}{\end{equation}}

\usepackage[obeyDraft]{todonotes}

\newcommand{\bal}{\begin{aligned}}
\newcommand{\enbal}{\end{aligned}}
\numberwithin{equation}{section} 
\newcommand{\cF}{{\mathcal F}}
\newcommand{\cA}{{\mathcal A}}
\newcommand{\cB}{{\mathcal B}}
\newcommand{\farc}{\frac}

\begin{document}

\title{\LARGE \textbf{Mean Field Games without Rational Expectations}\footnote{We thank Charles Bertucci, Larry Christiano, Jean-Michel Lasry and Peter Tankov for useful comments.}}
	\author{Benjamin Moll \thanks{Department of Economics, London School of Economics; London WC2A 2AE, UK; b.moll@lse.ac.uk. BM was partially supported by the Leverhulme Trust.}
		\and Lenya Ryzhik \thanks{Department of Mathematics, Stanford University; Stanford, CA, 94305,
		USA; ryzhik@stanford.edu. LR was partially supported by NSF grant DMS-2205497 and by ONR grant N00014-22-
1-2174.}
	}
\date{First version: March 2025	\\
	This version: January 2026 \\
\href{https://benjaminmoll.com/MFGRatX/}{[latest version]}}
\maketitle

\begin{abstract}
Mean Field Game (MFG) models implicitly assume ``rational expectations", meaning that the 
heterogeneous agents being modeled correctly know all relevant transition probabilities for the 
complex system they inhabit. When there is common noise, it becomes necessary to solve the 
``Master equation",   in which 
the infinite-dimensional density of agents is a state variable. The rational expectations assumption and the implication that agents solve Master equations is unrealistic in many applications. We show how to instead formulate MFGs with non-rational expectations. Departing from rational expectations is particularly relevant in ``MFGs with a low-dimensional coupling", i.e. MFGs in which agents' running reward function depends on the density only through low-dimensional functionals of this density. This happens, for example, in most macroeconomics MFGs in which these low-dimensional functionals have the interpretation of ``equilibrium prices." In MFGs with a low-dimensional coupling, departing from rational expectations allows for completely sidestepping the Master equation and for instead solving much simpler finite-dimensional HJB equations. We introduce an adaptive learning model as a particular example of non-rational expectations and discuss its properties.
\end{abstract}

\section{Introduction}

Economists and mathematicians cast models with a large number of interacting agents as Mean Field Games (MFGs), a coupled system of a backward-in-time Hamilton-Jacobi-Bellman (HJB) equation for agents' value function and a forward-in-time Fokker-Planck type equation for the agents' density. These equations describe the Nash equilibrium of a game played by a large number of agents experiencing fluctuations that are independent from each other. When there is common noise, the backward-forward stochastic coupling becomes more complicated and, to find their optimal strategy, the model agents need to solve a Master equation, that is, an equation in which the infinite-dimensional density is a state variable. This Master equation therefore suffers from an extreme version of the curse of dimensionality and has been nicknamed ``Monster equation" due to its complexity.\footnote{{Of course, this complexity does not arise in MFGs in which individual states take only a small number of possible values (say two or three) so that the Master equation is finite- and low-dimensional.}} {While this complexity makes the Master equation a fascinating mathematical object, it also limits its practical applicability.}

An underappreciated fact is that MFGs not only impose Nash equilibrium but also assume  ``rational expectations", meaning that the heterogeneous agents being modeled correctly know all relevant transition probabilities governing the complex system they inhabit. We argue that this assumption is unrealistic in many applications and show how to instead formulate a more general class of MFGs with non-rational expectations. Furthermore, we show that departing from rational expectations can drastically simplify the complexity of agents' optimization problems in certain applications and may allow sidestepping the unrealistically complex infinite-dimensional Master equation altogether.

After spelling out general MFGs without rational expectations, we focus on a class of MFGs we term ``MFGs with a low-dimensional coupling." In these MFGs, agents' running reward function depends on the density only through low-dimensional functionals of this density, for example a low-dimensional vector of moments of this density. That is, model agents do not directly ``care about" the density (in the sense that their rewards do not depend on it) and instead care only about the low-dimensional functionals.  Most MFGs in macroeconomics are examples of this class of MFGs, with these low-dimensional functionals corresponding to ``equilibrium prices."\footnote{In economics, these MFGs with coupling via low-dimensional equilibrium prices are known as ``heterogeneous agent models."} We show that, in MFGs with a low-dimensional coupling, departing from rational expectations generally results in a much simpler finite-dimensional HJB equation in place of the infinite-dimensional Master equation.

One of this paper's arguments is that MFGs with rational expectations and the Master equation are unrealistically complex as models of human decision making. MFGs with a low-dimensional coupling illustrate this clearly: under the rational expectations assumption, the low-dimensional coupling neither simplifies the formulation of the MFG nor that of the Master equation in any straightforward way. Indeed, to compute the low-dimensional functionals (``prices" in macroeconomics applications) in the rational expectations regime, one needs to compute the full density of the agents and there is no closed system that includes the ``prices" alone and not the full agents density. The agents being modeled are assumed to perform the same computations. It is for this reason that the Master equation is an infinite-dimensional PDE despite model agents only ``caring about" much lower-dimensional ``prices." The present paper criticizes the use of the Master equation in MFGs with a low-dimensional coupling and calls for developing alternative low-dimensional approximations that take advantage of these models' special structure. This part of the paper is a ``mathematics translation" using the language of partial differential equations of an economics paper \citep{moll-challenge} which criticizes the use of the rational expectations assumption in macroeconomics MFGs.\footnote{In macroeconomics MFGs, forward-looking  decision makers are assumed to forecast equilibrium prices by forecasting functionals of infinite-dimensional densities. But it seems self-evident that real-world households and firms do not forecast prices in this way and instead solve simpler  approximate problems.} We should also be clear that the main focus of this paper is on 
modeling and not on the mathematical analysis of the proposed models.

Our main results are contained in Section~\ref{sec:MFG_RE} which formulates MFGs without rational expectations. The preceding Sections~\ref{sec:MFG} to \ref{sec:nonRatX_simple} contain background material and building blocks that are useful for understanding such MFGs. In particular, before considering the full MFG case, Section  \ref{sec:nonRatX_simple} introduces the idea of departing from rational expectations in the simpler case of a single agent solving a stochastic control problem.

As we explain in Section~\ref{sec:MFG_RE}, the key feature of MFGs without rational expectations is that each agent uses a perceived trajectory of the future empirical density of the other agents that does not necessarily coincide with the density's actual equilibrium trajectory. This results in a model that, formally, still has the backward-forward feature of the MFG but, for a specified perceived trajectory of the density, at any given time, the agents' strategy can be computed solely by using the backward-in-time HJB equation, without resorting to the forward in time density equation. This system is the analogue of what economists call a ``temporary equilibrium" \citep{hicks-value,lindahl,grandmont,grandmont-palgrave}. After spelling out this more general system we show that, as expected, we recover the familiar backward-forward MFG system in the special case of rational expectations, i.e. when agents' perceived trajectory of the density of the other agents coincides with density's actual equilibrium trajectory. Analogously, we show how to formulate MFGs without rational expectations in the case with common noise and that we recover the Master equation in the special case with rational expectations.

In Section~\ref{sec:nonRE_lowD} we consider non-rational expectations in MFGs with a low-dimensional coupling. Each agent now uses a perceived trajectory of the future vector of ``prices" that does not necessarily coincide with the actual trajectory of equilibrium prices. The result is again a system in which, for a specified perceived trajectory of these prices, at any given time, the agents' strategy can be computed from the HJB equation without resorting to the Fokker-Planck equation. Importantly, in the case with common noise, for a specified perceived law of motion of future prices that imposes the Markov property, agents solve a simple finite-dimensional HJB equation. That is,  departing from rational expectations can completely sidestep the infinite-dimensional Master equation. 

Sections~\ref{sec:MFG_RE} and \ref{sec:nonRE_lowD} considered MFGs in which agents hold beliefs about future evolution or future prices that are specified outside the model (the ``temporary equilibrium" idea). In Section~\ref{sec:learning} we instead explain how such beliefs may be determined ``inside the model" via some form of learning. Specifically, we introduce an adaptive learning model that has the same key property as the models with exogenously-specified beliefs we just discussed: at any given time, and given the current prices, the agents strategy can be computed solely from the backward-in-time HJB equation, without resorting to the forward in time density equation. However, in contrast to the models with exogenous beliefs, with adaptive learning, agents update their beliefs in the face of new information so that, over time, their perceived trajectory of equilibrium prices may approximate the corresponding actual trajectory. Finally, we also discuss some other promising directions, in particular reinforcement learning and other stochastic approximation algorithms.

In Section~\ref{sec:discrete_time}, for the sake of exposition and completeness, we explain how the arguments of this paper can be adapted to the discrete-time case. Section \ref{sec:conclusion} concludes.

\section{Background: Mean Field Games and Rational Expectations\label{sec:MFG}}

In this section, we briefly review the basics of MFGs. \citet{lasry-lions}, \citet{cardaliaguet}, \citet{ryzhik}, 
\citet{carmona-delarue1,carmona-delarue2}, \citet{cardal-porr} 
and \citet{cardaliaguet-delarue-lasry-lions}  provide more 
complete treatments. 
We use the standard formulation in the MFG literature with small modifications explained below. 
We then briefly explain the rational expectations assumption.

\subsection{Backward-forward MFG system}\label{sec:back-for}

Let us recall the setup of mean-field games. Consider a system of $N\gg 1$ individual 
agents (players) at positions (states) $X_{i,t}\in\Rm^n$, $i=1,\dots,N$, 
with $0 \leq t \leq T$, where $T$ is a fixed terminal time that is sometimes taken as $T=+\infty$.

Given $t\ge 0$ and initial state $x\in\Rm^n$ (which differs across agents), 
each agent's state evolves according to the stochastic differential equation (SDE)
\begin{equation}\label{eq:SDE}
dX_{i,s} = \alpha_{i,s} ds + \sqrt{2 \nu} dB_{i,s}, \quad X_{i,t} = x,~~t\le s\le T.
\end{equation}
Here $\alpha_{i,s} \in A \subset \mathbb{R}^n$ is a control which is optimally chosen by each agent (in a way prescribed below), $B_{i,s}$ is a standard $n$-dimensional Brownian motion, and $\nu \geq 0$ is a parameter measuring the strength of the fluctuations affecting individual agents. 
The individual Brownian motions~$B_{i,s}$,~$i=1,\dots,N$, 
are independent and capture \emph{idiosyncratic} risk. The intuition is that this risk averages to a deterministic effective mean-field effect when $N\gg 1$, as far as the evolution of the agents' density is concerned.

To choose the control $\alpha_{i,s}$ in (\ref{eq:SDE}), each agent solves an optimization problem 
for the value function~$u_{N}: \mathbb{R}^n \times [0,T] \rightarrow \mathbb{R}$ defined by 
\begin{equation}\label{eq:objective-N}
u_{N}(x,t) = \max_{\alpha_i \in A} 
\mathbb{E} \left[\int_t^T e^{-\rho t} R(X_{i,s},\alpha_{i,s},m_{N}(s,\cdot))ds + 
e^{-\rho (T-t)} V(X_{i,T},m_{N}(T,\cdot)) \right] 
\end{equation}
subject to $X_{i,s}$ solving \eqref{eq:SDE}. Here,  
\be
m_N(t,x)=\sum_{i=1}^N\delta(x-X_{i,t}),~~x\in\Rm^n,
\ee
is the empirical measure of the collection of agents,  
$\rho \geq 0$ is a discount rate, 
$R(x,\alpha,m)$ is a running reward function that depends on 
the state $X_{i,s}$, the control $\alpha_{i,s}$ and the empirical
density $m_{N}(s,\cdot)$,  
and~$V(x,m)$ is a prescribed terminal value at time $T$. The expectation in (\ref{eq:objective-N}) is taken with respect to the idiosyncratic noises $B_{i,s}$, $1\le i\le N$.

In this setting, the dynamics of all agents are identical except 
that they experience different realizations of~$B_{i,t}$. That is, they are ex-ante identical (the functions $R,V,$ and so on are identical for all agents) but ex-post heterogeneous in~$X_{i,t}$ because of 
the different realizations of idiosyncratic risk $B_{i,t}$. One small modification to the 
standard formulation in the MFG literature is that here the 
agents maximize their objectives rather than minimize them.

Note that each agent's running reward $R$ depends on the overall system's state, 
the empirical density $m_{N}$. 
Agents optimally choose the control $\alpha$ taking the future evolution of $m_{N}$ as given. 
We denote the optimal policies, that is, the optimally chosen $\alpha$ in
(\ref{eq:objective-N}), by $\pi$ (see~(\ref{eq:policy_def}) below).  
Because~$R$ depends 
on $m_{N}$, so do the optimal policies $\pi$ and the value function $u_{N}$. In turn, 
the evolution of the empirical density $m_N$ depends on each agent's optimal policy
that appears in (\ref{eq:SDE}). We are thus considering the Nash equilibrium of a game between a large number of statistically identical players.

The backward-forward MFG system arises in the limit $N\to+\infty$ of a large number
of agents and is a coupled system of a Hamilton-Jacobi-Bellman equation for 
the limit $u(x,t)$ of the value functions $u_N(x,t)$ and a Fokker-Planck equation for 
the limiting empirical density of the agents~$m(x,t)$:
\begin{equation}\label{eq:BF}
\begin{split}
\rho u -  \partial_t u  &= H(x,\nabla u,m) + \nu \Delta u, \qquad \qquad \qquad \ \ \ \mbox{in} \ \mathbb{R}^n \times (0,T), \\
\partial_t m &= -\textrm{div}(\nabla_\lambda H(x,\nabla u,m) \ m) + \nu \Delta m, \quad \mbox{in} \ \mathbb{R}^n \times (0,T),\\
m(0)&=m_0, \quad u(x,T) = V(x,m(T)), \qquad \quad \ \mbox{in} \ \mathbb{R}^n.
\end{split}
\end{equation}
Here, $H$ is the Hamiltonian\footnote{We use $\lambda$ instead of the more standard $p$ to 
denote the dual "momentum" variable because $p$ will denote the price vector below.}
\begin{equation}\label{eq:hamiltonian}
H(x,\lambda,m) := \max_{\alpha \in A} \left\{R(x,\alpha,m) + \lambda \cdot \alpha \right\},~~\lambda\in\Rm^n, 
\end{equation}
and the policy function (optimal control) of each agent, defined as
\begin{equation}\label{eq:policy_def}
\pi(x,t) \equiv \arg \max_{\alpha \in A} \left\{ R(x,\alpha,m) + \alpha\cdot\nabla u(x,t) \right\},
\end{equation}
is given by
\begin{equation}\label{eq:control_hamiltonian}
\pi(x,t) = \nabla_\lambda H(x,\nabla u(x,t),m(t)).
\end{equation}

\subsection{Master equation with common noise}

We next introduce \emph{common noise} (in MFG terminology) or \emph{aggregate uncertainty} 
(in macroeconomics terminology). Typically, as, for example, in
\cite{cardaliaguet-delarue-lasry-lions,ahuja-2016,ahuja-2019}, this is done by directly adding 
an additional noise to the dynamics of individual agents in (\ref{eq:SDE})
\begin{equation}\label{eq:SDE_noise_MFG}
dX_{i,t} = \alpha_{i,t} dt + \sqrt{2 \nu} dB_{i,t} + \sqrt{2 \beta} dW_t, \quad X_{i,0} = x.
\end{equation}
Here, the extra noise $W_t$ is identical for all agents. In this  formulation, the agents density $m_t$ satisfies a \emph{stochastic} Fokker-Planck equation 
\be\label{25may1304}
d m_t = \left[-\textrm{div}(\nabla_\lambda H(x,\nabla_x u_t,m_t) \ m_t) + (\nu + \beta) \Delta_x m_t \right]dt - \textrm{div}(m_t \sqrt{2 \beta} dW_t),
\ee
rather than the standard Fokker-Planck equation that appears in (\ref{eq:BF}).
 
For the sake of simplicity of the presentation, we will discuss now a slightly different model, where the common noise does not directly affect the 
evolution of the individual agents themselves. Instead, in this model, the common noise directly affects the running reward
function. This will affect the optimal policy $\pi$ and thus the evolution of the
individual agents as well. \citet{bertucci-meynard1,bertucci-meynard2} use a similar approach.
 
Specifically, we introduce an additional state variable $Z_t \in \mathbb{R}^k$ 
(``the aggregate state"),
with some $k\ll N$, that evolves according to an~SDE:
\begin{equation}\label{25apr2710}
\begin{split}
dZ_t &= \mu_z(Z_t)dt+\sqrt{2\beta} dW_t, \quad Z_0 = z, 
\end{split}
\end{equation}
with some drift $\mu_z(z)$ and $\beta \geq 0$.
Here, $W_t$ is a standard $k$-dimensional \emph{common} Brownian motion that -- in contrast to the idiosyncratic Brownian motions $B_{i,t}$ -- affects \emph{all agents} simultaneously via $Z_t$. 
 

The main assumption is that the running reward function~$R$ 
in~\eqref{eq:objective-N} now depends on the aggregate state $Z_t$:
\begin{equation}\label{25apr2802}
u_{N}(x,t) = \max_{\alpha_i \in A} 
\mathbb{E} \left[\int_t^T e^{-\rho t} R(X_{i,s},Z_s,\alpha_{i,s},m_N(s,\cdot))ds + 
e^{-\rho (T-t)} V(X_{i,T},Z_T,m_N(T,\cdot)) \right]. 
\end{equation}
As usual in the common noise setting, the admissible controls $\alpha_{i,s}$ in
(\ref{25apr2802}), need to be $\cF_s$-measurable: they can not depend on the future.
This, of course, is also true for the optimal control $\pi(X_{i,t},Z_t,m_t,t)$.

If the reward function $R$ is non-separable between $Z_t$ and $\alpha_{i,t}$ (which is the relevant assumption in macroeconomics applications), 
the optimal policy $\pi$ will depend on the aggregate state $Z_t$ and therefore so will the dynamics of $X_{i,t}$:
\begin{equation}\label{eq:SDE_policy}
\bal
dX_{i,t}& = \pi(X_{i,t},Z_t,m_t,t) dt + \sqrt{2 \nu} dB_{i,t},\\ 
dZ_t &= \mu_z(Z_t)dt+\sqrt{2\beta} dW_t. 
\enbal
\end{equation}
Let us comment that this setting is a special case of the standard common noise formulation in~(\ref{eq:SDE_noise_MFG}) but with a degenerate diffusion.
To see this, introduce $\widetilde X_{i,t} \in \mathbb{R}^d$, with~$d=n+k$, with the two components 
\[
\widetilde X_{i,t} = \left[\begin{matrix} X_{i,t}\\ Z_{t} \end{matrix}\right].
\]
Note that the second component is identical for all agents. Then \eqref{eq:SDE_policy} is the special case of~(\ref{eq:SDE_noise_MFG})
in which the first component $X_{i,t}$ is only affected by the idiosyncratic noise 
but not by the common noise, and the second component $Z_t$ is only affected by the common noise but not by the idiosyncratic noise. Furthermore, the second
component is not controlled: $\alpha_{Z,t} \equiv 0$.

The state of the system (``the economy") is now a pair $(m_t,Z_t)$ which evolves as:
\begin{align}
\label{eq:dm}d m_t &= \left[-\textrm{div}_x(\nabla_\lambda H(x,Z_t,\nabla_x u_t,m_t) \ m_t) + \nu \Delta_x m_t \right]dt,\\
\label{eq:dZ}dZ_t &= \mu_z(Z_t)dt + \sqrt{2\beta} dW_t.
\end{align}
In contrast to the standard MFG formulation with a common noise \eqref{25may1304}, the Fokker-Planck equation (\ref{eq:dm}) for $m_t$ is not a stochastic
partial differential equation (SPDE) but a partial differential equation (PDE). Nevertheless, the solution to the system (\ref{eq:dm})-(\ref{eq:dZ}) is, formally,
an infinite-dimensional (degenerate) diffusion, and  $m_t$ itself is a stochastic object as it depends on the diffusion $Z_t$. Because $m_t$ is stochastic, it is now necessary to include it as a state variable in the agents' value function
$$U(x,z,m,t).$$
Note that $m \in \mathcal{P}(\mathbb{R}^n)$, the space of probability measures on with support in $\mathbb{R}^n$, which is an infinite-dimensional space.
Hence, the Master equation is a Hamilton-Jacobi-Bellman equation for the value function $U$ set in infinite-dimensional space:
\begin{equation}\label{eq:master}
\begin{split}
\rho U - \partial_t U =
& H(x,z,\nabla_x U,m) + \nu \Delta_x U + \beta \Delta_z U\\
& + \int_{\mathbb{R}^n} [\nabla_m U](y)  \underbrace{\left[-\textrm{div}_y(\nabla_\lambda H(y,z,\nabla_y U,m) \ m) + \nu \Delta_y m \right](y)}_{\textrm{drift of probability measure $m$ at point $y$ from \eqref{eq:dm}}} dm(y)\\
& \qquad \qquad \quad \ \mbox{in} \ \mathbb{R}^n \times \mathbb{R}^k \times \mathcal{P}(\mathbb{R}^n) \times (0,T)\\
U(x,z,m,T) =& V(x,z,m)\quad \mbox{in} \ \mathbb{R}^n \times \mathbb{R}^k \times \mathcal{P}(\mathbb{R}^n)
\end{split}
\end{equation}
Here $\nabla_m U$ denotes the derivative of $U$ with respect to the measure $m$ -- see \citet{cardaliaguet-delarue-lasry-lions} for a precise definition -- and $[\nabla_m U](y)$ denotes the derivative of $U$ with respect to $m$ \emph{at point $y$.} Also note that we set the drift of the aggregate state $\mu_z(z) \equiv 0$ for notational simplicity and we will continue to do so going forward.

The stochastic Fokker-Planck equation (\ref{25may1304}) that comes from the standard common noise MFG formulation (\ref{eq:SDE_noise_MFG}) 
further complicates the Master equation \eqref{eq:master} with additional second-order derivatives in $m$ \citep{cardaliaguet-delarue-lasry-lions}. The formulation in this section results in a simpler first-order Master equation and nests all typical macroeconomics applications.

\subsection{Rational expectations}
Rational expectations is a modeling assumption introduced by \citet{muth} and popularized in the 1970s by Bob Lucas, Ed Prescott, Tom Sargent and others. 
It has since been the standard assumption for modeling expectations in macroeconomics. See \citet{moll-challenge} for key references and a brief historical discussion.
The macroeconomics definition of rational expectations is as follows:
\emph{Agents have rational expectations if they form expectations over outcomes using the correct objective probability distributions of those outcomes. 
Hence, subjective probability distributions equal objective probability distributions.}

The rational expectations assumption is best thought of as a consistency requirement between expectations and model reality. Arguably a better name for the assumption is ``model-consistent expectations" \citep{simon-rationality}. 
Related, it is important to emphasize that ``rational expectations" is distinct from the concept of ``rationality" which, in economics terminology, simply means that agents maximize some objective function. While the MFGs we consider below relax rational expectations, all of them retain rationality.

In the context of MFGs, the rational expectations assumption is about the expectation operator $\mathbb{E}$ that appears in the 
optimization procedure for the objective functions \eqref{eq:objective-N} and \eqref{25apr2802}. If this expectation operator uses objective, model-consistent probability distributions 
for the behavior of all other agents as well as the idiosyncratic and common noise, then the MFG assumes rational expectations.

One key takeaway is that, as we explain in more detail in the next sections, \emph{all} existing Mean Field Games models in the mathematics literature (that we are aware of) implicitly assume rational expectations.

\section{Mean Field Games with a Low-dimensional Coupling}\label{sec:special_structure}

In Section~\ref{sec:MFG}, the running reward function $R(x,z,\alpha,m)$ 
that appears in the optimization problem~(\ref{25apr2802})
depends on the empirical measure $m(t)$ in a general, unrestricted fashion. 
However, in many applications, in particular in macroeconomics, 
this dependence is simpler:  the running reward function depends on the empirical
measure $m(t)$ only through a \emph{low-dimensional} vector $p_t \in \mathbb{R}^\ell$, with some fixed $\ell\ll N$, that is a functional of $m_t$.

That is, the running reward and the terminal condition in \eqref{25apr2802} are given by $\widetilde{R}(x,z,\alpha,p)$ and~$\widetilde{V}(x,z,p)$, 
so that agents optimize
\begin{equation}\label{objective_MFG_lowD}
u_{N}(x,t) = \max_{\alpha_i \in A} 
\mathbb{E} \left[\int_t^T e^{-\rho t} \widetilde R(X_{i,s},Z_s,\alpha_{i,s},p_s)ds + 
e^{-\rho (T-t)} \widetilde V(X_{i,T},Z_T,p_T) \right], 
\end{equation}
subject to $X_{i,s}$ solving \eqref{eq:SDE} and where
\be\label{25apr2704}
p_t= P^*(m_t,Z_t),
\ee
for a fixed functional 
\be
P^*: \mathcal{P}(\mathbb{R}^n) \times \mathbb{R}^k \rightarrow \mathbb{R}^\ell.
\ee
We refer to such MFGs as \emph{MFGs with a low-dimensional coupling}. Note that in MFGs with a low-dimensional coupling, model agents do not 
directly ``care about" the infinite-dimensional density $m_t$ in the sense that it does not enter their running reward functions or terminal conditions. 
Instead, they only ``care about" the much lower-dimensional vector $p_t$.

The next two subsections present 
applications from macroeconomics that take this form. Other examples are dynamic games with a large number 
of players in which the running reward depends only on particular moments of the distribution, say, the first 
moments 
\[
\overline{X}_j(m) = \int x_j dm(x,t), \quad j=1,...,n.
\]

Of course, an MFG with a low-dimensional coupling is just a special case of the general MFGs 
discussed in Section~\ref{sec:MFG} with the running reward and terminal value of the form 
\be
R(x,z,\alpha,m) = \widetilde R(x,z,\alpha,P^*(m,z)), \qquad V(x,z,m) = \widetilde V(x,z,P^*(m,z)).
\ee 
Therefore, under the rational expectations assumption, this low-dimensional coupling does not really simplify the analysis and the backward-forward MFG and Master 
equation still take the same form. However, as we will show below, low-dimensional coupling can drastically simplify the model's complexity when agents have non-rational 
expectations.

\subsection{Macroeconomics MFGs: low-dimensional coupling through prices}\label{sec:macroMFG}
Typical macroeconomics MFGs, known as ``heterogeneous agent models," are MFGs with a low-dimensional coupling. Usually, the running reward function $R(x,z,\alpha,m)$ depends 
on the empirical measure $m(t)$ only through a \emph{low-dimensional} price vector $p_t \in \mathbb{R}^\ell$, with some fixed $\ell\ll N$. 
The prices may represent the actual prices of goods, or correspond to wages or interest rates, that is, prices of other variables like labor and capital.
The underlying macroeconomic assumption is that the system
stays in what is known as a \emph{competitive equilibrium}. In that case, 
the prices are  set 
by the intersection of demand and supply (``market clearing") and are determined
by the empirical measure of the agents via a set of $\ell$ ``market clearing" conditions (demand equals supply for each of $\ell$ goods):
\be\label{25apr2702}
M(p_t, m_t,Z_t)=0,
\ee
with a given relation $M: \mathbb{R}^\ell \times \mathcal{P}(\mathbb{R}^n) \times \mathbb{R}^k \rightarrow \mathbb{R}^\ell$. Under the assumption that 
(\ref{25apr2702}) can be inverted, this gives rise to a unique mapping which 
takes the form in \eqref{25apr2704} and which we will refer to as the \emph{equilibrium price function}.

To summarize, the reward function depends on the measure~$m_t$ only 
through the (low-dimensional) 
price vector~$p_t$ that determines the optimal strategy in (\ref{eq:objective-N}).
In a competitive equilibrium, the price vector is directly related to $m_t$ either explicitly by (\ref{25apr2704}) or implicitly by~(\ref{25apr2702}).

\subsection{Example of a macroeconomics MFG}\label{sec:macro_example}
A typical example of a macroeconomics MFG that has this structure 
is the model described in \citet[][Section 5]{PDE-macro} which is a continuous-time version of 
the \citet{krusell-smith} model. The corresponding model without common noise is described in  \citet[][Section 2]{PDE-macro} and \citet{AHLLM} which is a continuous-time version of the Aiyagari-Bewley-Huggett model \citep{aiyagari,bewley,huggett}.

In this model, the state
of the agents is  characterized by their income and wealth, so that the agents' positions 
are parametrized by $x = (x_1,x_2) \in \mathbb{R}^2$. 
Here, $x_1$ is wealth and $x_2$ is income, and   $n=2$ in the setting of Section~\ref{sec:MFG}. 
Furthermore, there are $\ell=2$ prices~$p_t = (p_{1,t},p_{2,t}) \in \mathbb{R}^2$ where $p_{1,t}$ is the interest rate and $p_{2,t}$ the wage. There is $k=1$ aggregate state $Z_t \in \mathbb{R}$ which has the interpretation of the (logarithm of) productivity of a so-called ``representative firm" (see below).

Wealth and income evolve according to a system of SDE 
\be\label{eq:x1x2}
\bal
&dX_{1,i,t}   = (P_1^*(m_t,Z_t) X_{1,i,t} +P_2^*(m_t,Z_t)X_{2,i,t} - C_{i,t})dt, \\
&dX_{2,i,t}  = \mu(X_{2,i,t})dt+(2\nu)^{1/2}dB_{i,t}.
\enbal
\ee
Here, $C_{i,t}$ is the $i$-th agent's consumption that serves as a control in this setting
and $\mu$ is a drift coefficient. 
The price functionals~$P_1^*$ and~$P_2^*$ are the (scalar) interest rate and wage which depend on the measure $m$ via an equilibrium condition explained below. 
Agents choose their consumption to maximize a utility function 
\be\label{25apr2706}
\mathbb{E}\int_0^T e^{-\rho t} U(C_{i,t})dt \quad \mbox{where} \quad U(c) = \frac{c^{1-\gamma}}{1-\gamma}, \quad \gamma>0
\ee
subject to \eqref{eq:x1x2} and a state constraint $X_{1,i,t} \geq 0$.

The Hamiltonian  \eqref{eq:hamiltonian} is then 
\[
H(x,z,\lambda,m) = \max_{c} \left\{u(c) + \lambda_1(P_1^*(m,z) x_1 + P_2^*(m,z)x_2 - c) 
+ \lambda_2 \mu(x_2) \right\}.
\]
The Hamiltonian $H$ is non-linear and non-separable between $x,z,\lambda$ and $m$. At the same time,  it depends on $m$ only through the two-dimensional prices $P^*(m,z):\mathcal{P}(\mathbb{R}^2) \times \mathbb{R}\rightarrow \mathbb{R}^2$.

The price functionals (equilibrium wage and interest rate) are given in this model by
$$P_1^*(m,z) = \partial_{\overline{X}_1} F(\overline{X}_1(m),\overline{X}_2(m),z), \qquad P_2^*(m,z) = \partial_{\overline{X}_2} F(\overline{X}_1(m),\overline{X}_2(m),z)$$
where 
\begin{equation}\label{eq:production_function}
F(\overline{X}_1,\overline{X}_2,z) = e^z\sqrt{\overline{X}_1 \overline{X}_2}, \qquad \overline{X}_1(m) = \int x_1 dm(x,t), \qquad \overline{X}_2(m) = \int x_2 dm(x,t).
\end{equation}
The function $F$ has the interpretation of the production function of a so-called 
``representative" firm, $z$ that of (the logarithm of) the firm's productivity, and the derivatives $\partial_{\overline{X}_1} F$ 
and $\partial_{\overline{X}_2} F$ those of ``marginal products". The dependence of the prices merely on the first moments 
of $m$, that is, on~$\overline{X}_1(m)$ and~$\overline{X}_2(m)$,  is special to this particular application. Other macroeconomics applications feature (considerably) more complicated price functionals.

\subsection{The case without common noise}
For future references it will also be useful to spell out MFGs with a low-dimensional coupling but without common noise. This is simply the case in which neither the running reward $R$ nor the low-dimensional functional $p$ depend on the aggregate state $z$, i.e. this functional is simply given by $p_t= P^*(m_t)$ with $P^*: \mathcal{P}(\mathbb{R}^n) \rightarrow \mathbb{R}^\ell$. Equivalently, the case without common noise sets $Z_t=0$ for all $t$. In the macroeconomics MFG of Sections \ref{sec:macroMFG} and \ref{sec:macro_example}, the underlying assumption is that neither the running reward $R$ nor the market clearing condition $M$ depend on the aggregate state $Z$. The particular example in Section \ref{sec:macro_example} is the Aiyagari-Bewley-Huggett model analyzed in \citet{PDE-macro,AHLLM}.

\subsection{Backward-forward system and Master equation for MFGs with low-dimensional coupling}\label{sec:master_lowD}
The corresponding backward-forward MFG system in the case of a low-dimensional coupling takes exactly the same form \eqref{eq:BF} but with $R(x,\alpha,m) = \widetilde R(x,\alpha,P^*(m))$ in the Hamiltonian \eqref{eq:hamiltonian} and terminal value $V(x,m) = \widetilde V(x,P^*(m))$. The same is true for the Master equation for the value function $U(x,z,m,t)$ which takes the form \eqref{eq:master} but with $R(x,z,\alpha,m) = \widetilde R(x,z,\alpha,P^*(m,z))$ in the Hamiltonian and terminal value $V(x,z,m) = \widetilde V(x,z,P^*(m,z))$.


Note that the special structure of MFGs with a low-dimensional coupling neither simplifies the backward-forward MFG system nor the Master equation in any straightforward way. In particular the infinite-dimensional measure $m \in \mathcal{P}(\mathbb{R}^n)$ is still a state variable in the Master equation. What is noteworthy about this is that the Master equation is an infinite-dimensional
PDE despite model agents only “caring about” much lower-dimensional “prices.” As we explain in the following, the root cause of this feature is the rational expectations assumption.

\section{Non-rational Expectations in Simple Control Problems}\label{sec:nonRatX_simple}

Before specifying what rational expectations -- and departures from such expectations -- mean in the context of large systems of heterogeneous agents (i.e. MFGs) let us first consider the simpler case of a single agent solving a stochastic control problem. We turn to  MFGs in Section~\ref{sec:MFG_RE}.

\subsection{A simple stochastic control problem in an evolving environment}

Consider a single agent solving a stochastic control problem in a prescribed time-dependent environment: 
\begin{equation}\label{eq:objective_simple}
u(x,t) = \max_{\alpha \in A} \mathbb{E} \Big[\int_{t}^T e^{-\rho (\tau-t)} 
R(X_{\tau},\alpha_{\tau},\beta_\tau)d\tau + e^{-\rho (T-t)} V(X_{T}) \Big], 
\end{equation}
subject to $X_{\tau}$ solving a stochastic differential equation
\begin{equation}\label{eq:LOM_simple}
dX_{\tau} = \alpha_{\tau} d\tau + \sqrt{2 \nu} dB_{\tau}, \quad X_{t} = x,~~t\le \tau\le T.
\end{equation}
Here, $\alpha_{\tau}$ is the control used on the time interval $t\le\tau\le T$, 
and $\beta_\tau$ represents a known time-dependent environment that the agent cannot control. 

For future purposes, it will be useful to write the HJB equation for the value function $u$ in terms of the infinitesimal generator which summarizes the transition probabilities of the process for $X_t$:
\begin{equation}\label{25may1524}
\mathcal{A}_{\pi} := \pi\cdot\nabla + \nu \Delta,
\end{equation}
where $\pi(x,t)$ is the agent's policy. The HJB equation is then
\begin{equation}\label{eq:HJB_simple}
\rho u -  \partial_t u  = \max_{\alpha \in A} \left\{R(x,\alpha,\beta_t) + \mathcal{A}_\alpha u \right\}, \qquad \mbox{in} \ \mathbb{R}^n \times (0,T),
\end{equation}
with the terminal condition 
\be
u(x,T) = V(x), \qquad \hbox{in $\mathbb{R}^n$},
\ee
The  associated optimal policy function is 
\be\label{25may1502}
\pi(x,t) \equiv \arg \max_{\alpha \in A} \left\{ R(x,\alpha,\beta_t) + \alpha\cdot\nabla u(x,t) \right\}.
\ee

To explain what rational expectations mean in this setting, 
it is useful to spell out what this optimal control problem and the associated HJB equation look like without imposing the rational expectations assumption. 
In this model, expectations may be non-rational  in two respects: 
first, the agent may have incorrect beliefs about its own dynamics
for a given control~$\alpha_t$. Second, she may have incorrect beliefs about the environment $\beta_t$ in the future. We consider these two cases separately below.  

\subsection{Expectations about the evolution of the agent's state}\label{sec:beliefs_agent}

Let us fix a terminal time $T$ and for the moment fix some time $t\in(0,T)$. 
For simplicity of notation, we will  consider in this sub-section a time-independent environment, so that
\be
\beta_\tau \equiv\hbox{const} \quad \mbox{all} \ \tau.
\ee
As just discussed, the agent has rational expectations about the process for $X_{\tau}$ for $\tau>t$ if she forms expectations about $X_{\tau}$ using the correct objective transition probabilities. These transition probabilities are summarized via the infinitesimal generator. Rational expectations means that the agent's beliefs are summarized by the correct generator $\mathcal{A}_\pi$ defined in \eqref{25may1524} that determines the evolution of the actual state $X_t$ in accordance with (\ref{eq:LOM_simple}).

Non-rational expectations means that the agent has some other subjective beliefs about the future evolution of $X_{\tau}$ for $\tau>t$ summarized by a different generator $\widehat{\mathcal{A}}_{\pi}$. For example, the agent may believe that the state follows an alternative diffusion process
\begin{align}
\label{25may1402_simple}
d \widehat X_{\tau} &= \widehat{\mu}(\widehat{X}_{\tau},\alpha_{\tau},\tau)d\tau
+ \sqrt{2 \widehat{\nu}(\widehat{X}_{\tau},\alpha_{\tau},\tau)}dB_\tau,~~\tau\ge t,\\
\label{25may1404_simple}
\widehat X_{t}&=x
\end{align}
instead of the process~\eqref{eq:LOM_simple} in which case the generator is
\begin{equation}\label{eq:generator_nonRE}
\widehat{\mathcal{A}}_{\pi} := \widehat{\mu}(x,\pi,t)\cdot\nabla + \widehat{\nu}(x,\pi,t) \Delta.
\end{equation}
In particular, it may be the case that 
either~$\widehat{\mu}(x,\alpha,\tau) \neq \alpha$ 
or  $\widehat{\nu}(x,\alpha,\tau) \neq \nu$. This happens, for example, if the agent does not have
the full information about the idiosyncratic noise, or presumes the existence of an additional drift in the problem.

Agents' policies~$\pi(x,t)$ are now determined from the optimization problem
\begin{align}
\label{eq:expected_value_nonRE}
\widehat{u}(x,t) &= \max_{\alpha \in A} \widehat{\mathbb{E}} 
\Big[\int_{t}^T e^{-\rho (\tau-t)} R(\widehat{X}_{\tau},\alpha_{\tau})d\tau 
+ e^{-\rho (T-t)} V(\widehat{X}_{T}) \Big],
\end{align}
supplemented by (\ref{25may1402_simple})-(\ref{25may1404_simple}) 
for the evolution of $\widehat X_{\tau}$ on 
the same time interval~$t\le \tau\le T$. We denote the expectations operator by $\widehat{\mathbb{E}}$ to highlight that these subjective expectations are, 
in general, different from the objective (rational) expectations operator $\mathbb{E}$. That is, $\Em$ refers to the expectation
with respect to the trajectories generated by the evolution (\ref{eq:LOM_simple}), 
and~$\widehat\Em$ to those generated by the perceived evolution  (\ref{25may1402_simple}).

The HJB equation with non-rational expectations is therefore
\begin{align}
\label{eq:HJB_nonRE_simple}
\rho \widehat{u} - \partial_s \widehat{u}  &= \max_{\alpha \in A} \left\{ R(x,\alpha) 
+ \widehat{\mathcal{A}}_{\alpha}\widehat{u} \right\},  \quad \mbox{in} \ \mathbb{R}^n \times  (t,T)\\
\widehat{u}(x,T) &= V(x),\quad \mbox{in} \ \mathbb{R}^n. 
\end{align}
The associated policy with non-rational expectations is 
\be\label{actual_policy_simple}
\pi(x,t) = \arg \max_{\alpha \in A} 
\left\{ R(x,\alpha) + \alpha\cdot\nabla_x \widehat{u}(x,t) \right\},
\ee
rather than by (\ref{25may1502}). Note the difference between the two value functions that appear in (\ref{25may1502}) and (\ref{actual_policy_simple}).

%

Let us comment that the actual probability density $\rho(x,t)$ of the agent who follows the strategy $\pi(x,t)$ defined by (\ref{actual_policy_simple})
is
\be\label{25may1525}
\partial_t\rho=\cA_\pi^*\rho.
\ee
Note the difference between the operator $\widehat\cA_\pi$ that appears in the Hamilton-Jacobi-Bellman equation (\ref{eq:HJB_nonRE_simple}) 
and the operator $\cA_\pi$ whose adjoint appears in (\ref{25may1525}). 

\paragraph{The special case of rational expectations of the evolution of $X_t$.}
Rational expectations mean that the  generator~$\widehat{\cA}_\pi$ that appears in the HJB equation \eqref{eq:HJB_nonRE_simple} 
coincides with the actual generator~$\cA_\pi$. In other words, the perceived evolution (\ref{25may1402_simple}) coincides with the actual evolution (\ref{eq:LOM_simple}) for any given control $\alpha_{\tau,t}$. With this assumption, the HJB equation \eqref{eq:HJB_nonRE_simple} becomes \begin{equation}\label{eq:HJB_simple_operator}
\begin{split}
\rho u -  \partial_t u  &= \max_{\alpha \in A} \left\{ R(x,\alpha) + \mathcal{A}_{\alpha}{u} \right\}, \qquad      \mbox{in} \ \mathbb{R}^n \times (0,T), \\
u(x,T) &= V(x), \qquad \quad \ \mbox{in} \ \mathbb{R}^n.
\end{split}
\end{equation}
which is the same as (\ref{eq:HJB_simple}) for the case of constant $\beta_t$. Hence, for the special case of rational expectations, we have recovered the standard HJB equation.

\subsection{Expectations about the evolution of the environment}\label{sec:expectations_environment}

Let us now consider the same optimization problem but reintroduce the time-dependent environment $\beta_t$. 
At a time $t\in(0,T)$ the agent  has access to the past trajectory 
\be
\beta_{\le t}=\{\beta(t'):~0\le t'\le t\}.
\ee
She uses this information to make a prediction
\be\label{eq:perceived_LOM_beta}
\widehat\beta_{s,t}=\Theta(s,t;\beta_{\le t}),\qquad s>t,
\ee
with some given function $\Theta$ that depends on the running time $s$, the starting time $t$ and the path $\beta_{\le t}$ that was observed prior to the
time $t$. Non-rational expectations in this context mean that 
\be\label{25may1508}
\widehat\beta_{s,t} \neq \beta_s,\qquad\hbox{for some $s>t$},
\ee
i.e. that the agent's perceived trajectory of her external environment does not coincide with the environment's actual trajectory.

Let us assume for simplicity of notation that, while the environment may be predicted incorrectly, as in (\ref{25may1508}),
the perceived law of motion of the state $X_t$ is correct and the agent assumes that her trajectory is given by (\ref{eq:LOM_simple}):
\begin{equation}\label{25may1523}
dX_{\tau,t} = \alpha_{\tau,t} d\tau + \sqrt{2 \nu} dB_{\tau}, \quad X_{t,t} = x,~~t\le \tau\le T.
\end{equation}
In other words, we have both $\widehat\mu=\alpha$ and $\widehat\nu=\nu$ in (\ref{25may1402_simple}). 
Then, the agent is solving the optimization problem
\begin{equation}\label{25may1504}
u(x,t) = \max_{\alpha \in A} \mathbb{E} \Big[\int_{t}^T e^{-\rho (\tau-t)} 
R(X_{\tau,t},\alpha_{\tau,t},\widehat\beta_{\tau,t})d\tau + e^{-\rho (T-t)} V(X_{T,t}) \Big]
\end{equation}
subject to \eqref{25may1523} and \eqref{eq:perceived_LOM_beta}.\footnote{The value $u(x,t)$ also depends on the entire prediction $\widehat\beta_{\tau,t},\tau \geq t$ but we suppress this dependence for simplicity.} As the environment $\widehat\beta_{\tau,t}$ now depends both on the running time $\tau$ and on the  starting time~$t$, in order to
formulate the corresponding Hamilton-Jacobi-Bellman equation, it is convenient to fix~$t\in(0,T)$ and introduce an auxiliary optimization problem
\begin{equation}\label{25may1506}
\widehat u(x,s;t) = \max_{\alpha \in A} \mathbb{E} \Big[\int_{s}^T e^{-\rho (\tau-s)} 
R(X_{\tau,s},\alpha_{\tau,s},\widehat\beta_{\tau,t})d\tau + e^{-\rho (T-s)} V(X_{T,t}) \Big],~~t<s<T,
\end{equation}
subject to \eqref{25may1523} and \eqref{eq:perceived_LOM_beta}. Note that the prediction of the future environment $\widehat\beta_{\tau,t}$ is made at the time $t$ and does not depend on the intermediate times $s$.
Then, the perceived future optimal policy of the agent is, similarly to (\ref{25may1502}), given by
\be\label{25may1518}
\widehat{\pi}(x,s;t) = \arg \max_{\alpha \in A} \left\{ R(x,\alpha,\widehat\beta_{s,t}) + \alpha\cdot\nabla_x \widehat{u}(x,s;t) \right\},~~t\le s\le T.
\ee
In particular, at the time $s=t$ we have 
\be 
\label{25may1516}
\pi(x,t) = \arg \max_{\alpha \in A} \left\{ R(x,\alpha,\widehat\beta_{t,t}) 
+ \alpha\cdot\nabla_x \widehat{u}(x,t;t) \right\}= 
\arg \max_{\alpha \in A} \left\{ R(x,\alpha,\beta_t) + \alpha\cdot\nabla_x {u}(x,t) \right\},
\ee
with 
\be\label{25may1520}
u(x,t)=\widehat u(x,t;t),
\ee
being the optimal value defined by (\ref{25may1504}).
This is the actual policy that the agents are following. As a result of this optimization problem, an agent's state evolves as a diffusion
\be\label{25may1512}
dX_{t} = \pi(X_{t},t)dt + \sqrt{2\nu}dB_{t},
\ee
with the policy $\pi(x,t)$ given by (\ref{25may1516}).

To find the function $u(x,t)$ one needs to solve a backward-in-time  Hamilton-Jacobi-Bellman equation for the value function $\widehat u(x,s;t)$.
It takes the form 
\begin{equation}\label{25may1514}
\bal
\rho \widehat{u}(x,s;t) - \partial_s \widehat{u}(x,s;t)  &= \max_{\alpha \in A} \left\{ R(x,\alpha,\widehat\beta_{s,t}) 
+ \mathcal{A}_{\alpha}\widehat{u}(x,s;t) \right\},  \quad \mbox{in} \ \mathbb{R}^n \times  (t,T), \\
\widehat{u}(x,T;t) &= V(x),\quad \mbox{in} \ \mathbb{R}^n .
\enbal
\end{equation}
Once the function $\widehat u(x,s;t)$ is found, one defines $u(x,t)$ by (\ref{25may1520}), which, in turn, gives the policy
(\ref{25may1516}) that appears in the actual dynamics (\ref{25may1512}).

\paragraph{The special case of rational expectations for the evolution of the environment.}
Rational expectations on the evolution of the environment mean that 
\be\label{25may1521}
\widehat\beta_{s,t}\equiv \beta_s,\qquad\hbox{for all $s>t$.}
\ee
Note that, for deterministic variables like the one considered here, rational expectations simply means that agents have \emph{perfect foresight} about the evolution of 
these variables. If that is the case, then the function $\widehat u(x,s;t)$ that solves (\ref{25may1514})
does not depend on the starting time $t$ and neither does the policy $\pi(x,s;t)$ that also
appears in (\ref{25may1514}). 
Under the assumption~(\ref{25may1521}), the HJB equation \eqref{25may1514} therefore becomes the standard  HJB equation (\ref{eq:HJB_simple}) with 
\be
\widehat u(x,s;t)=u(x,s).
\ee

\section{Mean Field Games without Rational Expectations}\label{sec:MFG_RE}

As already noted, \emph{all} existing Mean Field Games models in the mathematics literature we are aware of implicitly assume rational expectations. We now explain this in more detail and show how to formulate Mean Field Games without rational expectations.  We first consider the case without common noise, that is, the backward-forward MFG  system~\eqref{eq:BF}, and then cover the Master equation \eqref{eq:master}.

\subsection{Non-rational expectations in the backward-forward MFG system}\label{sec:nonRE-back-for}

We now show how to formulate a generalization of the backward-forward MFG system of 
Section \ref{sec:back-for} without making the rational expectations assumption. We focus on the 
more interesting case of non-rational expectations about the evolution of the agent's external environment, 
in this case the evolution of the density $m$. The case of non-rational expectations about agents' own 
states $X_{i,t}$ is analogous to the treatment in Section \ref{sec:beliefs_agent}. As in the preceding 
section, after spelling out the model with general (not necessarily rational) expectations, we show that we 
recover the standard backward-forward MFG system in the special case with rational expectations.

Let us fix a terminal time $T$ and for the moment fix some time $t\in(0,T)$.
We assume that to predict the future empirical density of the other
agents, an individual agent uses a perceived law of motion 
\begin{align}
\label{25may1402}
\partial_s \widehat m(x,s;t) &= \mathcal{B}^* \widehat m(x,s;t),~~s\ge t,\\
\label{25may1404}
\widehat m(x,t;t)&=m(x,t)
\end{align}
Here, $\widehat m(x,s;t)$  is a prediction for the empirical measure 
of other agents for $s\ge t$, and 
$\cB$ is the generator of a Markov process that an agent believes the other agents are following. 
The initial condition in (\ref{25may1404}) at $s=t$ comes
from the actual observed density $m(x,t)$ at the time $t$. Note that the agents are constantly updating their 
prediction $\widehat m(x,s;t)$  for a given future time $s$, by changing the
initial condition in~(\ref{25may1404}) as~$t$ grows. In other words, $\widehat m(x,s;t)$ and $\widehat m(x,s;t')$ are, in general, different if~$s>t>t'$. 
Going forward, we sometimes write $\widehat{m}_{s,t}$ for conciseness.

Note that this setup is exactly analogous to the case of non-rational expectations about the evolution of the environment in Section \ref{sec:expectations_environment}. Therefore, so is the remainder of the discussion. Similarly to \eqref{25may1506}, agents' policies~$\pi(x,t)$ are determined from the perceived 
optimization problem
\begin{align}
\label{eq:expected_valueXXX}
\widehat{u}(x,s;t) &= \max_{\alpha_i \in A} \mathbb{E} 
\Big[\int_{s}^T e^{-\rho (\tau-s)} R(X_{i,\tau},\alpha_{\tau,s},\widehat{m}_{\tau,t})d\tau + 
e^{-\rho (T-s)} V(X_{i,T,s},\widehat{m}_{T,t}) \Big] ,~~t\le s\le T,
\end{align}
subject to \eqref{eq:SDE} for the evolution of $X_{i,t}$ and (\ref{25may1402})-(\ref{25may1404}) for the evolution of $\widehat m_{s,t}$ on the time interval $t\le s\le T$.
As a result, the perceived future optimal policy of the agent is, similarly to (\ref{25may1518}), given by
\be\label{25may1410}
\widehat{\pi}(x,s;t) = \arg \max_{\alpha \in A} \left\{ R(x,\alpha,\widehat{m}_{s,t}) 
+ \alpha\cdot\nabla_x \widehat{u}(x,s;t) \right\},~~t\le s\le T.
\ee
At the time $s=t$ we have 
\be 
\label{eq:policy_learningXXX}
\bal\pi(x,t) &=\widehat{\pi}(x,t;t)= 
\arg \max_{\alpha \in A} \left\{ R(x,\alpha,m_t) + \alpha\cdot\nabla_x \widetilde{u}(x,t) \right\}.
\enbal
\ee
with $\widetilde u(x,t)=\widehat u(x,t;t)$. This is the actual policy that the agents are following. 

As a result of this optimization problem, an idiosyncratic state evolves as a diffusion
\be\label{25may1416}
dX_{i,t} = \pi(X_{i,t},t)dt + \sqrt{2\nu}dB_{i,t}.
\ee
To summarize, agents believe that for all $s>t$ the distribution of the other agents will evolve according to
(\ref{25may1402})-(\ref{25may1404}),  and this is what shows up in the continuation values for $s\ge t$ in the definition~(\ref{eq:expected_valueXXX})
of the perceived value function $\widehat u(x,s;t)$. 
However, when they choose their actual policy $\pi(x,t)$ at time $t$, given by  (\ref{eq:policy_learningXXX}), 
they use the actual realized density $m_t$ (which will generally differ from any previous estimates of $\widehat m_{t,t'}$ with $t'<t$)
and the perceived value function $\widetilde u(x,t)$.

Agents' actual policies $\pi(x,t)$ defined in (\ref{eq:policy_learningXXX}) give rise to the generator 
\begin{equation}\label{eq:generatorXXX}
\mathcal{A}_{\pi} := \pi\cdot\nabla + \nu \Delta.
\end{equation}
that determines the evolution of the actual density $m(x,t)$  in accordance with (\ref{25may1416}): 
\be\label{25may1408}
\partial_tm=\cA_\pi^* m,
\ee
which is different from the evolution (\ref{25may1402})-(\ref{25may1404}) for the perceived future density unless $\cB=\cA_\pi$.

Now, the backward-foward system of equations for the value function $\widehat u(s,x;t)$,
the perceived density $\widehat m(x,s;t)$ and the actual density $m(x,t)$  becomes:
\begin{equation}\label{eq:BFMFG_nonrational}
\bal
\rho \widehat{u}(x,s;t) - \partial_s \widehat{u}(x,s;t)  &= \max_{\alpha \in A} 
\left\{R(x,\alpha, \widehat m_{s,t}) 
+ \mathcal{A}_{\alpha}\widehat{u}(x,s;t)\right\}, \quad \qquad \quad ~\mbox{in} \ \mathbb{R}^n \times  (t,T), \\
\partial_s \widehat{m}(x,s;t) &= \mathcal{B}^* \widehat{m}(x,s;t), \ \ \ \quad \qquad \qquad \qquad \qquad \qquad \qquad \mbox{in} \ \mathbb{R}^n \times (t,T),\\
\pi(x,t) &= \arg \max_{\alpha \in A} \left\{ R(x,\alpha,m(t)) 
+ \alpha\cdot\nabla_x \widehat{u}(x,t;t) \right\}, \\
\partial_t m &= \mathcal{A}^*_\pi m, \ \ \qquad \qquad \qquad \qquad \qquad \qquad \qquad \qquad  \mbox{in} \ \mathbb{R}^n \times (0,T),\\
m(x,0)&=m_0(x), \ \ \widehat m(x,t;t) = m(x,t),\ \ \widehat{u}(x,T;t) = V(x,\widehat{m}_{T,t}) \ \ \mbox{in} \ \mathbb{R}^n.
\enbal
\end{equation}
The system \eqref{eq:BFMFG_nonrational} is the analogue of what economists call a ``temporary equilibrium".

\textbf{Definition:} \emph{Temporary equilibrium} at a particular time $t$ is defined as allocations 
and policies such that (i) agents optimize \emph{given expectations of future variables} 
(including the density) that are specified in the model but that are \emph{not necessarily rational}, 
(ii) the economy is in Nash equilibrium at time $t$.

This idea was originally developed contemporaneously by \citet{hicks-value} and \citet{lindahl}, and
has been further developed by \citet{grandmont,grandmont-palgrave}

\paragraph{Remark: connection to the Master equation.} 
Note that the system (\ref{eq:BFMFG_nonrational}) 
could also be written in terms of a Master equation, with the measure~$m$ as a state variable.
However, this would defeat the purpose of this approach: the advantage of (\ref{eq:BFMFG_nonrational})
is exactly that we do not have to compute the dynamics in the infinite-dimensional space of probability
measures but rather only find the solutions to  (\ref{eq:BFMFG_nonrational}) for the measures $m_t$ that
one  encounters in the course of actual evolution.

\paragraph{Rational expectations in the Backward-Forward MFG System.} Rational expectations in the context of the model we have discussed above 
mean that the  generator~$\cB$ that appears in the equations 
(\ref{25may1402}) and (\ref{eq:BFMFG_nonrational}) for the perceived density $\widehat m(x,s;t)$
coincides with the actual generator~$\cA_\pi$ that appears in the evolution
equation (\ref{eq:generatorXXX}) for the actual density~$m(x,t)$:
\be\label{25may1412}
\cB=\cA_\pi.
\ee
Let us now show that with this assumption, the system (\ref{eq:BFMFG_nonrational}) 
reduces to the familiar MFG system~(\ref{eq:BF}) that we write 
in the form 
\begin{equation}\label{eq:BFXXX}
\begin{split}
\rho u -  \partial_t u  &= \max_{\alpha \in A} \left\{ R(x,\alpha,m(t)) 
+ \mathcal{A}_{\alpha}{u}\right\}, \qquad      \mbox{in} \ \mathbb{R}^n \times (0,T), \\
\pi(x,t) &= \arg \max_{\alpha \in A} \left\{ R(x,\alpha,m(t)) + \alpha\cdot\nabla u(x,t) \right\}, \\
\partial_t m &= \mathcal{A}_{{\pi}}^*m,  \quad \qquad \qquad \qquad \qquad \qquad 
~~~~\mbox{in} \ \mathbb{R}^n \times (0,T),\\
m(0)&=m_0, \quad u(x,T) = V(x,m(T)), \quad \quad \ \mbox{in} \ \mathbb{R}^n.
\end{split}
\end{equation}
Indeed, if~(\ref{25may1412}) holds, we deduce from the second and fourth
equations in (\ref{eq:BFMFG_nonrational}), together with the 
initial condition $\widehat m(x,t;t)=m(x,t)$ that
\be\label{25may1418}
\widehat m(x,s;t)=m(x,s),~~\hbox{for all $t\in[0,T]$ and $s\in[t,T]$,}
\ee
so that the perceived future density $\widehat m(x,s;t)$ does not depend on the 
starting time $t$ and coincides with the actual density $m(x,s)$.  
The perceived value function $\widehat u(x,s;t)$ is then independent of~$t$ as well, so that we have
\be\label{25may1420}
\widehat u(x,s;t)=\widehat u(x,s;s)=\widetilde u(x,s).
\ee 
The perceived policies $\widehat\pi(x,s;t)$ are also independent of $t$ because
we can use (\ref{25may1418}) and (\ref{25may1420}) to write 
\be
\bal
\widehat{\pi}(x,s;t) &= \arg \max_{\alpha \in A} 
\left\{ R(x,\alpha,\widehat{m}_{s,t}) + \alpha\cdot\nabla_x \widehat{u}(x,s;t) \right\}\\
&=
\arg \max_{\alpha \in A} \left\{ R(x,\alpha,m(s)) + \alpha\cdot\nabla_x \widehat{u}(x,s;s) \right\}
=\pi(x,s).
\enbal
\ee
It follows from the above and the first equation in (\ref{eq:BFMFG_nonrational})
that the value function
$\widetilde u(s,x)$ satisfies the backward HJB equation
\be
\bal
\rho \widetilde{u}(x,s) - \partial_s \widetilde{u}(x,s)  &
= R(x,\pi(x,s),m(s)) 
+ \mathcal{A}_{{\pi}}\widetilde{u}(x,s)  \quad \mbox{in} \ \mathbb{R}^n \times  (s,T). 
\enbal
\ee
Together with the forward Fokker-Planck equation 
\be
\partial_t m = \mathcal{A}^*_\pi m \qquad  \mbox{in} \ \mathbb{R}^n \times (0,T)
\ee
we have recovered the MFG system \eqref{eq:BFXXX} or, equivalently, (\ref{eq:BF}).
 
This shows clearly that the backward-forward MFG system \eqref{eq:BF} implicitly assumes that 
agents have rational expectations 
about the process $X_{i,t}$ for all agents~$i=1,...,N$, and hence for the evolution of the density $m(x,t)$.
 
\paragraph{Remark: non-rational expectations about agents' own states.} As already noted, one could also allow for non-rational expectations about agents' own states $X_{i,t}$. Analogous to the treatment in Section \ref{sec:beliefs_agent}, this simply involves replacing the generator $\mathcal{A}_{\alpha}$ in the first equation in~\eqref{eq:BFMFG_nonrational} by some other generator $\widehat\mathcal{A}_{\alpha}$ capturing each agent's beliefs about the evolution of her own state, i.e. a perceived law of motion like \eqref{25may1402_simple} and \eqref{25may1404_simple}.

\paragraph{Remark: heterogeneous beliefs.}
It is reasonable to expect that different agents may have different beliefs about the evolution of the density of the other agents. In fact, substantial belief heterogeneity is one of the most prevalent empirical findings in the macroeconomics literature on household and firm expectations, see the discussion and references in \citet{moll-challenge} (the finding is often summarized under the terminology ``disagreement"). In our context, different individual agents may use different forms of the operator $\cB$, that appears in (\ref{25may1402}) to predict the future evolution of the density of the other agents. 
This will, in turn, affect the actual generator~$\cA_\pi$ governing the evolution of the actual density $m(t)$ of the agents. We will revisit this issue below in the discussion of adaptive learning in Section \ref{sec:learning_no_noise}.

\paragraph{Remark: the unrealism of rational expectations.} Rational expectations imposes that agents know the correct objective transition probabilities not only for their own individual states but also for the evolution of every other agents' states, i.e. for the entire complex system they inhabit as a whole. Specifically, the implicit assumption is that each agent \emph{knows all other agents' optimal policies} $\pi(y,\cdot)$ for all $y \in \mathbb{R}^n$ and \emph{uses these to (correctly) forecast the evolution of everyone else's state and hence the measure $m$!} For these reasons, the rational expectations assumption is arguably a stretch in complex environments like MFG models.

\paragraph{Remark: the evolution of the density with heterogeneous agents.}
It is also worth pointing out that rational expectations assumption (\ref{25may1412}) 
is the reason why the operator $\mathcal{A}^*_\pi$ in the Fokker-Planck equation in \eqref{eq:BFXXX} 
is the adjoint of the operator $\mathcal{A}_\pi$ in the HJB equation 
in \eqref{eq:BFXXX}. Without rational expectations this would not be the case: 
while the Fokker-Planck equation for the actual density $m(x,t)$ in (\ref{eq:BFMFG_nonrational})
necessarily features (the adjoint of) the correct generator~$\mathcal{A}_\pi$ because it 
reflects the actual realized dynamics of $X_{i,t}$ and the associated evolution of the measure $m$, 
the HJB equation for the perceived value function~$\widehat u(s,x;t)$ in
(\ref{eq:BFMFG_nonrational}) features $\mathcal{A}_\pi$ only under the rational expectations assumption (\ref{25may1412}).
Otherwise it is driven by the perceived strategy $\widehat\pi$.

The rational expectations assumption (\ref{25may1412}) may seem to superficially simplify a 
``complicated" system (\ref{eq:BFMFG_nonrational}) to a ``simpler" (and, definitely, shorter) system 
(\ref{eq:BFXXX}). However, while the latter may produce a higher value function for the agents, it comes
at a high computational cost. If the operator $\cB$ that appears in the evolution equation
for the perceived density $\widehat m(x,s;t)$ in~(\ref{eq:BFMFG_nonrational}) is independent
from the perceived value function $\widehat u(x,s;t)$ then (\ref{eq:BFMFG_nonrational}) may be solved
in a single pass for a fixed $t<T$.
First, one solves the forward in time equation
\be
\bal
\partial_s \widehat{m}(x,s;t)& = \mathcal{B}^* \widehat{m}(x,s;t),~~ 
\mbox{in} \ \mathbb{R}^n \times (t,T),\\
\widehat m(x,t;t)&=m(x,t),
\enbal
\ee
to find the perceived future density of the agents. 
This is followed by solving the backward-in-time Hamilton-Jacobi-Bellman equation 
\be
\bal
\rho \widehat{u}(x,s;t) - \partial_s \widehat{u}(x,s;t)  &= \max_{\alpha \in A} \left\{R(x,\alpha,\widehat m_{s,t}) 
+ \mathcal{A}_{\alpha}\widehat{u}(x,s;t)\right\},  \quad \mbox{in} \ \mathbb{R}^n \times  (t,T) ,\\
\widehat{u}(x,T;t)&= V(x,\widehat{m}(T,t)),
\enbal
\ee
to find the perceived value function $\widehat u(x,s;t)$.  
This requires no iterations that are normally used to solve the forward-backward MFG
problems. With $\widehat u(x,t;t)$ in hand, one can compute the actual policy $\pi(x,t)$ 
given by (\ref{eq:policy_learningXXX}),
and continue the evolution of the actual agents density $m(x,t)$ forward in time. 
Thus, from the computational-cost point of view, dropping the rational expectations
assumption may end up being beneficial.

\subsection{Non-rational expectations in Mean Field Games with common noise}\label{sec:nonRE_MFG_common}

We now repeat the exercise for the case with common noise and show how to formulate MFGs without rational expectations in this case. To this end, recall that the evolution of the states~$(X_{i,t},Z_t)$ is given by \eqref{eq:SDE_policy} where $\pi(X_{i,t},Z_t,m_t,t)$ is the optimal policy with common noise.

Analogous to \eqref{25may1402} and \eqref{25may1404}, with non-rational expectations, the agent incorrectly believes that the density and aggregate state evolve according to
\be\label{25may2602}
\bal
d \widehat{m}_{s,t} &= \mathcal{B}_{\widehat{Z}_s}^* \widehat{m}_{s,t} ds, \quad s \geq t\\
d \widehat Z_{s,t} &= \sqrt{2 \widehat{\beta}(\widehat{Z}_{s,t},\tau)}dW_s,~~s\ge t,\\
\enbal
\ee
with $\widehat{m}_{t,t}=m_t$ and $\widehat Z_{t,t}=Z_t$, instead of the correct evolution \eqref{eq:dm} and \eqref{eq:dZ}. As already noted, we set the drift $\mu_z(z) \equiv 0$ in  \eqref{eq:SDE_policy} for simplicity.
Note that, in general, the perceived generator $\mathcal{B}_{\widehat{Z}_s}$ may depend on the perceived aggregate state $\widehat{Z}_s$.

This leads to the following infinite-dimensional Master equation:
\begin{equation}\label{eq:master_nonRE}
\begin{split}
\rho \widehat{U} - \partial_t \widehat{U} &=
 \max_{\alpha \in A} \left\{ R(x,z,\alpha,m) + \mathcal{A}_{\alpha,z} \widehat{U} \right\} + \widehat{\beta}(z,t) \Delta_z \widehat{U} + \int_{\mathbb{R}^n} [\nabla_m \widehat{U}](y)  [\mathcal{B}^*_z m](y) dm(y)\\
& \qquad \qquad \qquad \ \ \mbox{in} \ \mathbb{R}^n \times \mathbb{R}^k \times \mathcal{P}(\mathbb{R}^n) \times (0,T)\\
\widehat{U}(x,z,m,T) &= V(x,z,m)\quad \mbox{in} \ \mathbb{R}^n \times \mathbb{R}^k \times \mathcal{P}(\mathbb{R}^n)
\end{split}
\end{equation}
As before, $\mathcal{A}_{\pi,z}$ in the term $\mathcal{A}_{\pi,z} \widehat U$ summarizes agents' beliefs about the evolution of the individual state $X_{i,t}$. 
In contrast, $[\mathcal{B}^*_z m](y)$ summarizes their beliefs about evolution of measure~$m$ at point $y$ which may, in general, differ from the actual evolution.

After solving for the optimal policies $\pi(x,z,m,t)$, the evolution of the actual density of the agents can be found from the coupled system of SDE
\begin{align}
\label{eq:dynamics_RE}d m_t = \mathcal{A}^*_{\pi,Z_t} m_t dt, \qquad dZ_t = \sqrt{2\beta} dW_t.
\end{align}
with
\[
\mathcal{A}_{\pi,z} = \pi(x,z,m,t)\nabla_x + \nu \Delta_x.
\]

\paragraph{Rational expectations.} Exactly as in Section \ref{sec:nonRE-back-for}, rational expectations mean that the  generator~$\cB$ that appears in the equations 
(\ref{25may2602}) and (\ref{eq:master_nonRE}) for the perceived density $\widehat m_{s,t}$ and the perceived value function $\widehat U$ coincides with the actual generator~$\cA_{\pi,z}$ that appears in the evolution
equation (\ref{eq:dynamics_RE}) for the actual density~$m(x,t)$:
\be
\cB_z=\cA_{\pi,z}, 
\ee
and, in addition, that the perceived diffusivity for the aggregate state $Z_t$ equals the actual diffusivity, $\widehat{\beta}(z,t)=\beta$.  In that case, \eqref{eq:master_nonRE} becomes 
\begin{equation}\label{eq:master_RE}
\begin{split}
\rho U - \partial_t U &=
\max_{\alpha \in A} \left\{ R(x,z,\alpha,m) + \mathcal{A}_{\alpha} U\right\} + \beta \Delta_z U + \int_{\mathbb{R}^n} [\nabla_m U](y)  [\mathcal{A}^*_{\pi,z} m](y) dm(y)\\
& \qquad \qquad \qquad \ \ \mbox{in} \ \mathbb{R}^n \times \mathbb{R}^k \times \mathcal{P}(\mathbb{R}^n) \times (0,T)\\
U(x,z,m,T) &= V(x,z,m)\quad \mbox{in} \ \mathbb{R}^n \times \mathbb{R}^k \times \mathcal{P}(\mathbb{R}^n)\\
\mathcal{A}_{\pi,z} &= \pi\nabla_x + \nu \Delta_x,\\
\pi(x,z,m,t) &= \arg \max_{\alpha \in A} \left\{ R(x,z,\alpha,m) + \mathcal{A}_{\alpha,z} U\right\}.
\end{split}
\end{equation}
Importantly, analogously to the discussion in the preceding subsection, the implicit assumption underlying rational expectations is that each agent knows not only the correct stochastic process for their own individual state $X_{i,t}$ but also that of the state of all other agents, i.e. their optimal policies in the future.

\paragraph{Remark: Nash equilibrium, ``common knowledge", and the computational complexity of the Master equation.} MFGs impose Nash equilibrium, i.e. that each agent plays a best response to their prediction of every other agents' strategy. With rational expectations, each agent \emph{knows} the other agents' strategies and then solves the infinite-dimensional Master equation to compute this best response. In fact, the Master equation not only imposes that all agents know each other's strategies, but they also know that they all know, and so on, ad infinitum. This assumption is called ``common knowledge" in the economics literature. In that sense, the policies coming from the solution to the Master equation are optimal for a given individual agent only if the other agents also have rational expectations and if everyone knows that they do. Without that social compact in which everyone predicts an identical future, solving the infinite-dimensional Master equation is suboptimal and may even be harmful. While requiring such ``common knowledge" is a common issue for Nash equilibria, the setting of the infinite-dimensional Master equation is somewhat special because computing the Nash equilibrium  requires a huge computational cost that may not be accessible to all agents in the system.

\section{Non-rational expectations in MFGs with low-dimensional coupling}\label{sec:nonRE_lowD}

In this section, we consider MFGs that have the special structure in Section \ref{sec:special_structure}: agents' rewards depend on the measure $m$ only through a low-dimensional functional.

\subsection{The case without common noise}

The case without common noise is analogous to Section \ref{sec:nonRE-back-for} and we therefore cover it only briefly. As in Section \ref{sec:special_structure}, 
the running reward and terminal value are given by $\widetilde{R}(X_{i,t},\alpha_{i,s},p_t)$ and~$\widetilde{V}(X_{i,t},p_t)$ for a low-dimensional vector 
$p_t \in \mathbb{R}^\ell$ with $p_t = P^*(m_t)$. Going forward we drop the tildes from $\widetilde{R}$ and $\widetilde{V}$ for notational simplicity.

With rational expectations, agents understand the dependence of $p_t$ on $m_t$ and therefore use the functional $P^*$ together with the correct evolution for $m_t$ to 
predict future values of~$p_t$. That is, agents have perfect foresight over the future trajectory of $p_t$. With non-rational expectations, 
agents instead perceive some other trajectory
$$
\widehat p_{s,t}=\Theta(s,t;p_{\le t}),\quad s>t, \quad \mbox{with} \ p_{t,t} = p(t).
$$
In particular note that agents' perceived trajectory of future ``prices" $\widehat p_{s,t},s>t$ generally depends on past realizations $p_{\le t}$ which model 
agents have already observed by the time $t$.

\paragraph{Benchmark with rational expectations.} We first spell out the backward-forward MFG system with rational expectations using the generator notation introduced in the preceding Section. We will refer back to this benchmark system at various future points of the paper. The system is:
\begin{equation}
\begin{split}\label{eq:BFMFG_foresight}
\rho u - \partial_t u  &= \max_{\alpha \in A} \left\{ R(x,\alpha,p(t)) + \mathcal{A}_{\alpha}u \right\} \qquad \mbox{in} \ \mathbb{R}^n \times (0,T) \\
\pi(x,t) &= \arg \max_{\alpha \in A} \left\{ R(x,\alpha,p(t)) + \alpha\nabla u(x,t) \right\},\\
\partial_t m &= \mathcal{A}^*_\pi m \quad \qquad \qquad \qquad \qquad \qquad \ \ \mbox{in} \ \mathbb{R}^n \times (0,T),\\
p(t) &= P^*(m(t)) \qquad \qquad \qquad \qquad \qquad \mbox{in} \ (0,T)\\
m(0)&=m_0, \quad u(x,T) =  V(x,P^*(m_T))\quad \mbox{in} \ \mathbb{R}^n. 
\end{split}
\end{equation}

\paragraph{Non-rational expectations.} Analogously to Section \ref{sec:nonRE-back-for}, the backward-forward MFG system without rational expectations is
\begin{equation}\label{eq:BFMFG_coupling_nonrational}
\bal
\rho \widehat{u}(x,s;t) - \partial_s \widehat{u}(x,s;t)  &= \max_{\alpha \in A} \left\{ R(x,\alpha,\widehat{p}_{s,t}) 
+ \mathcal{A}_{\alpha}\widehat{u}(x,s;t) \right\},  \quad \mbox{in} \ \mathbb{R}^n \times  (t,T), \\
\widehat p_{s,t} &=\Theta(s,t;p_{\le t}), \qquad \qquad \qquad \qquad \qquad ~~ \mbox{in} \ (t,T),\\
\pi(x,t) &= \arg \max_{\alpha \in A} \left\{ R(x,\alpha,p(t)) +
 \alpha\cdot\nabla_x \widehat{u}(x,t;t) \right\}, \\
\partial_t m &= \mathcal{A}^*_\pi m, \qquad \qquad \qquad \qquad \qquad \qquad ~~~~~  \mbox{in} \ \mathbb{R}^n \times (0,T),\\
p(t) &= P^*(m(t)), \ \qquad \qquad \qquad \qquad \qquad ~~~~~  \mbox{in} \ (0,T) \\
m(0,x)&=m_0(x),\quad \widehat{u}(x,T;t) = V(x,\widehat{p}_{T,t}),\qquad \mbox{in} \ \mathbb{R}^n .
\enbal
\end{equation}

\subsection{Common noise: sidestepping the Master equation in MFGs with a low-dimensional coupling}\label{sec:sidestepping_lowD}

In the presence of a common noise $Z_t$, as in Section \ref{sec:special_structure}, the running reward is $R(X_{i,t},Z_t,\alpha_{i,s},p_t)$ for a low-dimensional vector $p_t= P^*(m_t,Z_t)\in \mathbb{R}^\ell$. With rational expectations, agents understand the dependence of $p_t$ on $m_t$ and $Z_t$ and therefore use the functional $P^*$ together with the correct stochastic processes for $m_t$ and $Z_t$ to predict future values of $p_t$.
This leads to the Master equation (\ref{eq:master_RE_CE}), with essentially no simplifications despite the low-dimensional coupling.

With non-rational expectations agents instead perceive some other stochastic process for the pair~$(p_t,Z_t)$. 
In the simplest case, they simply perceive $p_t$ to evolve according to a completely exogenous stochastic process
\begin{equation}\label{eq:price_process_simple}
d\widehat{p}_{s,t} = \widehat{\mu}_p(\widehat{p}_{s,t})ds + \widehat{\sigma}_p(\widehat{p}_{s,t})dW_t, \ s \geq t \qquad \widehat{p}_{t,t} = p_t
\end{equation}
In more complicated cases, agents perceive a joint stochastic process for $p_t,Z_t$ and other variables (which could, in principle, include $m_t$).

In the case of agents perceiving the simple process \eqref{eq:price_process_simple}, instead of writing a Master equation, we can write a much simpler, 
standard finite-dimensional HJB equation for a value function $\widehat{U}(x,z,p,t)$. Denoting the generator corresponding to \eqref{eq:price_process_simple} by $\widehat \cA_p$ we have
\begin{equation}\label{eq:master_nonRE_CE}
\begin{split}
\rho \widehat U - \partial_t \widehat U &=
\max_{\alpha \in A} \left\{R(x,z,\alpha,p) + \mathcal{A}_\alpha \widehat U \right\} + \beta \Delta_z \widehat U + \widehat{\cA}_p \widehat{U} \quad  
\mbox{in} \ \mathbb{R}^n \times \mathbb{R}^k \times \mathbb{R}^\ell \times (0,T),\\
U(x,z,p,T) &=  V(x,z,p) \ \quad \qquad \qquad \qquad  \qquad \qquad   \qquad \qquad ~~  \mbox{in} \ \mathbb{R}^n \times \mathbb{R}^k \times \mathbb{R}^\ell,\\
\mathcal{A}_\alpha &= \alpha\nabla_x + \nu \Delta_x.
\end{split}
\end{equation}
Therefore, in MFGs with a low-dimensional coupling, departing from rational expectations can completely sidestep the Master equation. Of course, the case considered 
here is just an illustrative example. In particular, note that the perceived law of motion \eqref{eq:price_process_simple} is specified completely ``outside the model" 
which leaves open the question of where this perceived law of motion ``comes from" in the first place.

\subsection{The trouble with the Master equation in MFGs with a low-dimensional coupling \label{sec:trouble}}

The corresponding Master equation for the value function $U(x,z,m,t)$ is instead
\begin{equation}\label{eq:master_RE_CE}
\begin{split}
\rho U - \partial_t U &=
\max_{\alpha \in A} \left\{ R(x,z,\alpha,P^*(m,z)) + \mathcal{A}_\alpha U \right\} + \beta \Delta_z U + \int_{\mathbb{R}^n} [\nabla_m U](y)  [\mathcal{A}^*_{\pi,z} m](y) dm(y)\\
& \qquad \qquad \qquad \ \ \mbox{in} \ \mathbb{R}^n \times \mathbb{R}^k \times \mathcal{P}(\mathbb{R}^n) \times (0,T)\\
U(x,z,m,T) &= V(x,z,m)\quad \mbox{in} \ \mathbb{R}^n \times \mathbb{R}^k \times \mathcal{P}(\mathbb{R}^n)\\
\mathcal{A}_{\pi,z} &= \pi\nabla_x + \nu \Delta_x  \\
\pi(x,z,m,t) &= \arg \max_{\alpha \in A} \left\{ R(x,z,\alpha,P^*(m,z)) + 
\alpha\cdot\nabla_x U(x,z,m,t) \right\} 
\end{split}.
\end{equation}
This is the same Master equation as described in Section \ref{sec:master_lowD} but using the generator notation used in the present section. As already noted there,   
special structure of MFGs with a low-dimensional coupling does not simplify the Master equation in any straightforward way.
In particular the infinite-dimensional measure $m \in \mathcal{P}(\mathbb{R}^n)$ is still a state variable in agents' value function.

The reason this happens is the rational expectations assumption. Intuitively, because agents are forward-looking, they need to forecast future prices $p_t$, a low-
dimensional object. But they understand that $p_t$ depends on the infinite-dimensional measure $m_t$ via the functional~\eqref{25apr2704}. Therefore, agents 
forecast the measure $m_t$ in order to forecast prices $p_t$. Furthermore, as before, they forecast $m_t$ using their knowledge of all other agents' policies $\pi(y)$. Note that all of this happens despite agents not even directly ``caring about" the distribution $m_t$.

Related, note that actual equilibrium prices $p_t$ do not follow a Markov process (if they did, one could, write a finite-dimensional HJB equation with prices $p$ as the state variables also in the rational expectations case); instead only $(m_t,Z_t)$ has the Markov property and prices are instead a complicated non-linear functional of this Markov state. Agents with rational expectations therefore (unrealistically) forecast the Markov state $(m_t,Z_t)$ in order to forecast the non-Markovian prices $p_t$.

Considering the case of macroeconomics MFGs \citet{moll-challenge} argues that we should not make our lives so hard. It seems self-evident that real-world households and firms do not forecast prices by forecasting cross-sectional distributions and instead solve simpler problems. Instead of solving ``Monster equations" we should replace the rational expectations assumption and solve the simpler equations corresponding to households' and firms' actual price-forecasting behavior.

\section{A way forward: learning in MFGs}\label{sec:learning}

As we have argued, MFGs with rational expectations and the Master equation are unrealistically complex as models of human decision making. We have also seen that 
departing from rational expectations may hold the promise of sidestepping the infinite-dimensional Master equation altogether. However, in the MFGs with non-rational 
expectations we have considered thus far, agents held beliefs about future evolution or future prices that were specified outside the model (the ``temporary equilibrium" 
idea). Therefore, these specified beliefs may end up being completely disconnected from the actual evolution of these equilibrium objects, i.e. there may be a disconnect 
between beliefs and ``model reality" and agents' expectations may be systematically disappointed. 
A related issue is that a model with exogenously specified beliefs is subject to a version of the so-called `` Lucas critique" \citep{lucas-critique}: when there is a change in economic policy (which would typically correspond to a change in a model parameter), one should expect agents' beliefs to change as well and this belief updating needs to be modeled.
 
How can we model non-rational expectations that are endogenous to the actual equilibrium prices but that, nevertheless, sidestep the Master equation and allow for computing standard finite-dimensional HJB equations for agents' value functions? Put differently: how can we formulate, in a systematic way, models of agents' behavior in situations with a low-dimensional coupling that lead to equations that (i)~approximate agents' real-world behavior, and (ii) sidestep computing the solutions to a Master equation with the infinite-dimensional state $m \in \mathcal{P}(\mathbb{R}^n)$ and the associated curse of dimensionality? This is the challenge posed by \citet{moll-challenge}.

A natural answer is to add  some form of \emph{learning} to the model. That is, instead of imposing -- as the rational expectations assumption does -- 
that agents \emph{know} the correct (and extremely complex) transition probabilities of equilibrium prices, we instead impose that agents \emph{learn} about these 
transition probabilities over time. This approach has a long tradition in the economics literature, typically in the form 
of ``least-squares learning" \citep{bray,marcet-sargent,evans-honkapohja}, and has recently been applied in the MFG literature to the case without 
common noise \citep{lauriere-etal-RL,lauriere-etal-survey,xu-etal-RL,bertucci-incomplete}, mostly in the form of reinforcement learning.

Before proceeding to describing this way forward, let us also note that one promise of modeling 
learning in this way is to ``kill two birds with one stone'': to develop variants of MFGs with a low-dimensional coupling that are easier to compute and analyze, while, at the same time, making these models more realistic and more likely to generate interesting macroeconomic phenomena.

\subsection{Adaptive learning without common noise}\label{sec:learning_no_noise}

Adaptive learning is a special case of the non-rational beliefs about the future we have discussed in Sections~\ref{sec:MFG_RE} and~ \ref{sec:nonRE_lowD} above. 
More specifically, this is a special form of the future
prices predictor~$\Theta(s,t;p_{\le t})$ that appears in the system (\ref{eq:BFMFG_coupling_nonrational}). One simple version of adaptive learning is least-squares learning \citep{bray,marcet-sargent,evans-honkapohja}. \citet{jacobson} implements such an approach in a heterogeneous-agent model.

In the simplest version of least-squares learning, at any time $t$, agents simply expect prices to remain constant at some value $\bar{p}$, i.e. $p_s = \bar{p}$ for all $s > t$. However, they update their estimate of this constant value over time as they collect new data on actual realized prices and we denote agents' time-$t$ estimate of $\bar{p}$ by $\widehat{p}_t$. Specifically, agents compute $\widehat{p}_t$ as a simple backward-looking average:
$$
\widehat{p}_t = \frac{1}{t}\int_0^t p_s ds.
$$
Differentiating, we have 
\begin{equation}\label{eq:pdot_LS}
\dot{\widehat{p}}_t = \frac{1}{t}\left(p_t - \widehat{p}_t \right). 
\end{equation}
Other forms of learning besides least-squares learning are possible as well, such as, for example, an ordinary differential equation
of the form 
\be\label{25may3002}
\dot{\widehat{p}}_t = L(p_t,\widehat{p}_t).
\ee
For example, the ODE
\[
\dot{\widehat{p}}_t = \alpha\left(p_t - \widehat{p}_t \right),
\]
with a constant $\alpha>0$ rather than the factor of $1/t$, as in (\ref{eq:pdot_LS}), is what is called ``adaptive expectations" \citep{cagan}. 
Note that expected prices will generally differ from actual prices, i.e. there is no longer perfect foresight.

In more complicated versions of least-squares learning, agents have a parametric ``perceived law of motion" (PLM) of prices
\begin{equation}\label{eq:general_learning}
\dot{p}_t = \mu_p(p_t,\theta) 
\end{equation}
where $\theta\in \mathbb{R}^d$ is a parameter vector which parameterizes their beliefs about the evolution of~$p_t$. For example the PLM function $\mu_p$ could be linear: 
\[
\dot{p}_t = \theta_0 + \theta_1 p_t.
\]
Agents then update their estimate $\widehat{\theta}_t \in \mathbb{R}^d$ 
of the parameter vector $\theta$
recursively over time
\begin{equation}\label{eq:general_learning_theta}
\dot{\widehat{\theta}}_t = L(p_t,\widehat{\theta}_t), \quad \widehat{\theta}_0 \ \ \mbox{given},
\end{equation}
which plays a similar role to the ODE (\ref{25may3002}). For example, if $G$ is linear, $\widehat{\theta}_t$ could be a backward-looking least-squares estimator.
 Note that $p_t$ that appear in  (\ref{eq:general_learning_theta}) 
are the actual observed prices at the time $t$.

Agents' policies $\pi(x,t)$ are determined from the following optimization problem 
that is analogous to (\ref{eq:expected_valueXXX}). We fix $t\in[0,T]$ and for each $s\in[t,T]$ consider the perceived value function 
\begin{align}
\label{eq:expected_value}\widehat{u}(x,s;t) &= \max_{\alpha_i \in A} \mathbb{E} 
\left[\int_s^T e^{-\rho (\tau-s)} R(X_{i,\tau},\alpha_{i,\tau},\widehat p_{\tau;t})d\tau 
+ e^{-\rho (T-s)} V(X_{i,T},\widehat p_{T,t}) \right].
\end{align}
subject to \eqref{eq:SDE} for the evolution of $X_{i,t}$ and where $\widehat p_{\tau;t}$ are the perceived future prices that evolve according to  
\be\label{25may2804}
\bal
\farc{d\widehat p_{\tau;t}}{d\tau}& = \mu_p(\widehat p_{\tau;t},\widehat\theta_t), ~~t\le\tau\le T, \\
\widehat p_{t;t}&=p_t.
\enbal
\ee
Here, $p_t$ are the actual observed prices at the time $t$, and 
$\widehat\theta_t$ is the solution to (\ref{eq:general_learning_theta}). We emphasize that the parameter~$\widehat\theta_t$ in (\ref{25may2804})  is fixed for $t\le\tau\le T$ 
and does not depend on $\tau$. 
That is, from the point of view of an agent at time $t$, $\widehat\theta_t$ is just a fixed parameter vector that parameterizes 
the perceived law of motion $\mu_p$ of future prices $\widehat p_{\tau;t}$, 
with~$\tau\geq t$. Note that the dependence of the perceived value
function $\widehat u(x,s;t)$ on the time $t$ is solely through the parameter $\theta_t$ that appears in the perceived evolution  (\ref{25may2804})
of the future prices.

The optimization problem \eqref{eq:expected_value} gives rise to the perceived 
future policy
\be
\bal
\label{eq:policy_learning}
\widehat{\pi}(x,s;t) &= \arg \max_{\alpha \in A} \left\{ R(x,\alpha,\widehat p_{s;t}) 
+ \alpha\cdot\nabla_x \widehat{u}(x,s;t) \right\}.
\enbal
\ee
The actual policy that the agents are following is, on the other hand, 
\be\label{25may2806}
\pi(x,t)=\widehat\pi(x,t;t)= \arg \max_{\alpha \in A} \left\{ R(x,\alpha,p_t) 
+ \alpha\cdot\nabla_x \widehat{u}(x,t;t) \right\}.
\ee
As above, the interpretation is that agents believe that for all $s>t$ prices will evolve according to (\ref{25may2804}) and this is what shows up in their continuation values. However, when they choose their policy at time $t$, they see the actual realized prices $p_t$ (which will generally differ from any previous estimates of $p_t$).

Going forward we drop the hat from $\widehat \theta_t$ for notational simplicity but keep in mind that this is really a time-varying estimate of the parameter $\theta$ in \eqref{eq:general_learning}. With this notation in hand, the backward-foward MFG system with adaptive learning becomes
the following version of (\ref{eq:BFMFG_coupling_nonrational})
\begin{equation}\label{25may2810}
\bal
\rho \widehat{u}(x,s;t) - \partial_s \widehat{u}(x,s;t)  
&= \max_{\alpha \in A} \left\{ R(x,\alpha,\widehat p_{s;t}) + \mathcal{A}_{\alpha}\widehat{u}(x,s;t) \right\},\quad \mbox{in} \ \mathbb{R}^n \times  (t,T), \\
\farc{d\widehat p_{s;t}}{ds}& = \mu_p(\widehat p_{s;t},\theta_t), \quad \qquad \qquad \qquad \qquad \qquad \mbox{in} \ (t,T), \\
\widehat p_{t;t}&=p_t,\\
\pi(x,t) &= \arg \max_{\alpha \in A} \left\{ R(x,\alpha,p_t) +
 \alpha\cdot\nabla_x \widehat{u}(x,t;t) \right\}, \\
\partial_t m &= \mathcal{A}^*_\pi m, \ \ \ \qquad \qquad \qquad \qquad \qquad \qquad ~~~  
\mbox{in} \ \mathbb{R}^n \times (0,T),\\
\dot{\theta}_t &= L(p_t, \theta_t), \ \ \ \quad \theta_0 \ \ \mbox{given},
\qquad \qquad   ~~~~~~~~~  \mbox{in} \ (0,T),\\
p_t &= P^*(m(t)), \ \ \ \qquad \qquad \qquad \qquad \qquad ~~~~  \mbox{in} \ (0,T), \\
m(0,x)&=m_0(x),\quad \widehat{u}(x,T;t) = V(x,\widehat{p}_{T,t}), \ \qquad \mbox{in} \ \mathbb{R}^n .
\enbal
\end{equation}
Here, in the context of adaptive learning, the low dimensionality of the coupling is crucial:
instead of a highly complex infinite-dimensional analog of the Master equation, 
we get a finite-dimensional system (\ref{25may2810}).
The forward in time equation for $m(t,x)$ in the system (\ref{25may2810}) is coupled to the 
backward in time equation for the perceived value function
$\widehat u(x,s;t)$
solely via the actual price $p_t$ that appears in the strategy $\pi(x,t)$ and in the parameter $\theta_t$ that appears in the Hamilton-Jacobi-Bellman
equation for $\widehat{u}(x,s;t)$. This coupling is not as severe as in MFGs with rational expectations 
since, as we have mentioned previously, solution to (\ref{25may2810}) 
does not require any iterations. It does require, however, to solve, for each $t\in[0,T]$ a separate Hamilton-Jacobi-Bellman equation on time interval $t\le s\le T$.

\paragraph{The Hamilton-Jacobi-Bellman equation in the price space.}
An alternative way of writing \eqref{25may2810} is to introduce the perceived value 
function $\widehat U(x,p,s;\theta)$ defined for all prices $p \in \Rm^\ell$ and parameters $\theta \in \Rm^d$:
\begin{align}
\label{25jun108}\widehat{U}(x,p,s;\theta) &= \max_{\alpha_i \in A} \mathbb{E} 
\left[\int_s^T e^{-\rho (\tau-s)} R(X_{i,\tau,{s}},\alpha_{i,\tau,{s}},\widehat p_{{\tau},s})d\tau 
+ e^{-\rho (T-s)} V(X_{i,T,{s}},\widehat p_{T,s}) \right],
\end{align}
defined for $s\le T$. Note that the perceived value function no longer depends on 
the time $t$ at which the parameter $\theta_t$ is fixed for $s>t$.  As we have mentioned, the
only dependence in  (\ref{eq:expected_value}) on $t$ comes from the parameter $\theta_t$ that is fixed for $s>t$. Here, the agents are allowed to take various 
values of $\theta$, which, in turn, allows
us to get rid of the dependence on $t$.  The perceived future prices $\widehat p_{\tau,s}$ evolve according to a generalization of
(\ref{25may2804})
\be\label{25jun112}
\bal
\farc{d\widehat p_{\tau,s}}{d\tau}& = \mu_p(\widehat p_{\tau,s},\theta), ~~s\le\tau\le T, \\
\widehat p_{s,s}&=p,
\enbal
\ee
but now with the parameter $\theta$ set to the value that appears in the argument of $\widehat u(x,p,\theta,s)$ in~(\ref{25jun108}).

The perceived value function $\widehat U(x,p,s;\theta)$ defined in \eqref{25jun108} solves a single HJB equation 
\be\label{master_learning_no_noise}
\bal
\rho \widehat{U}(x,p,s;\theta) - \partial_s \widehat{U}(x,p,s;\theta)  &= H(x,p,\nabla_x \widehat{U} ) 
+ \nu \Delta_x \widehat{U}(x,p,s;\theta)+\mu_p(p,\theta)\cdot\nabla_p\widehat U(x,p,s;\theta), \\
H(x,p,\lambda) &= \max_{\alpha \in A} \left\{ R(x,\alpha,p) + 
\alpha\cdot \lambda \right\},
\enbal
\ee
with  $(x,p,s)\in \ \mathbb{R}^n \times\Rm^\ell \times  (0,T)$ and $\theta \in \Rm^d$. 

Note that the arguments of the value function include not only the low-dimensional prices $p \in \Rm^\ell$ but also the  "fairly high-dimensional" parameter vector $\theta \in \Rm^d$. While including $p$ comes at a low computational cost, including $\theta$ comes at a significantly higher additional computational cost. However, note that equations (\ref{master_learning_no_noise}) for different $\theta$ are decoupled from each other. As we will remark below, this is something one can exploit to reduce computational costs.

The HJB equation (\ref{master_learning_no_noise}) for the perceived value function gives rise to the perceived policy
\be
\bal
\label{25jun114}
\widehat{\pi}(x,p,s;\theta) &= \arg \max_{\alpha \in A} \left\{ R(x,\alpha,p) + \alpha\cdot\nabla_x \widehat{U}(x,p,s;\theta) \right\}.
\enbal
\ee
Analogously to above, the actual time-$t$ policy is then given by the perceived policy evaluated at the time-$t$ price $p_t$ and parameter $\theta_t$:
\be\label{25jun121_XXX}
\pi(x,t)=\widehat\pi(x,p_t,t;\theta_t).
\ee
The prices $p_t$ in the right side of (\ref{25jun121_XXX}) are given by
\be\label{25jun122_XXX}
p_t = P^*(m(t)).
\ee
The measure $m(x,t)$ and the parameter $\theta_t$, in turn, solve the forward-in-time problems
\be
\bal\label{25jun123_XXX}
\partial_t m(x,t) &= \mathcal{A}^*_\pi m(x,t),\\
\dot{\theta}_t &= L(p_t, \theta_t),
\enbal
\ee
with the policy $\pi$ given by (\ref{25jun121_XXX}) and with initial conditions $m_0(x)$ and $\theta_0$. The system (\ref{25jun121_XXX})-(\ref{25jun123_XXX}) 
is driven by the solution to (\ref{master_learning_no_noise})
via the policy $\widehat\pi$ that appears in (\ref{25jun121_XXX}) . 

\paragraph{Remark on computational cost of (\ref{master_learning_no_noise})}
The price for solving a single standard HJB equation for $\widehat{U}$ is the need to solve (\ref{master_learning_no_noise}) for all~$\theta$ in the region of interest. Thus, the computational complexity of this formulation is controlled by the dimension $d$ of the parameter space~$\theta \in \Rm^d$ which may be considerably higher than in \eqref{25may2810}. However, one can take advantage of the fact that the HJB equations~\eqref{master_learning_no_noise} are decoupled for different $\theta$ here. In particular, this means that one can solve~\eqref{master_learning_no_noise} only for the $\theta$'s one actually encounters. Specifically, starting from an initial condition $\theta_0$, solve the HJB equation \eqref{master_learning_no_noise} for $\theta=\theta_0$; then update $\theta$ according to \eqref{25jun123_XXX}, solve the HJB equation again for this new $\theta$, and so on.

\paragraph{Internalized learning}
In the HJB equation (\ref{master_learning_no_noise}), agents do not take into account that the parameter vector $\widehat\theta_t$ evolves according to the learning rule  \eqref{eq:general_learning_theta}. Hence, in (\ref{master_learning_no_noise}) learning is ``external" to agents. An alternative approach pursued in part the economics literature \citep[e.g.][]{christiano-eichenbaum-johannsen} is to instead assume that agents ``internalize" learning, taking into account the evolution of $\widehat\theta_t$. The corresponding variant of (\ref{master_learning_no_noise}) 
with internalized learning is:
\be\label{master_learning_no_noise_internal}
\bal
\rho \widehat{U} - \partial_s \widehat{U}  &= H(x,p,\nabla_x \widehat{U} ) 
+ \nu \Delta_x \widehat{U}+\mu_p(p,\theta)\cdot\nabla_p\widehat U  + L(p,\theta)\nabla_\theta \widehat U , \\
H(x,p,\lambda) &= \max_{\alpha \in A} \left\{ R(x,\alpha,p) + 
\alpha\cdot \lambda \right\},
\enbal
\ee
with  $(x,p,s)\in \ \mathbb{R}^n \times\Rm^\ell \times  (0,T)$ and $\theta \in \Rm^d$.
Both formulations are sensible from an economic perspective and simply correspond to different assumptions about agents' sophistication.

\paragraph{Belief heterogeneity.} In the system (\ref{25may2810}), the parameter 
$\theta_t$, that couples its forward and backward in time components, 
summarizes agents' beliefs about the evolution of the future prices. Although these beliefs 
vary over time, and~$\theta_t$ evolves in (\ref{25may2810}) according to a differential equation,
the assumption leading to (\ref{25may2810}) is that all agents share the \emph{same} beliefs $\theta_t$. As remarked in Section \ref{sec:nonRE-back-for}, it is natural to allow for heterogeneity in beliefs which is an important feature of real-world data on empirical measures of such beliefs.

A convenient feature of the formulation \eqref{master_learning_no_noise} in which $\theta$ is a state variable is that it is easy to extend to the case of belief heterogeneity which we model as follows.
Beliefs $\theta$ differ across the population and there is 
an initial joint density of states and beliefs $m_0(x,\theta)$. Starting from time $t=0$, beliefs evolve according to the following generalization of \eqref{eq:general_learning_theta} which allows for learning to depend on the individual state $X_{i,t}$ as well: 
\be
\dot{\theta}_{i,t} = L(p_t,X_{i,t},\theta_{i,t}), \quad \theta_{i,0} \sim m_0(x,\theta).
\ee
The value function of an agent with beliefs $\theta$ still satisfies the same HJB equation \eqref{master_learning_no_noise}. The actual policy of that agent is then given by 
\be\label{25jun116}
\pi(x,\theta,t)=\widehat\pi(x,p_t,t;\theta)= \arg \max_{\alpha \in A} \left\{ R(x,\alpha,p_t) + \alpha\cdot\nabla_x \widehat{u}(x,p_t,t;\theta) \right\},
\ee
meaning that their trajectories are 
\[
dX_{i,t} = \pi(x,\theta,t) d\tau + \sqrt{2 \nu} dB_{i,\tau}.
\]
The prices $p_t$ on the right side of (\ref{25jun116}) are still given by
\be\label{25jun122}
p_t = P^*(m(t)).
\ee
The measure $m(x,\theta,t)$, in turn, solves the forward in time problem
\be\label{25jun123}
\partial_t m(x,\theta,t) = \mathcal{A}^*_\pi m(x,\theta,t)-\hbox{div}_\theta(L(p_t,x,\theta)m(x,\theta,t)),
\ee
with the policy $\pi$ given by (\ref{25jun116}).

It is worth contrasting the system (\ref{25jun116})-(\ref{25jun123}) with belief heterogeneity with its counterpart (\ref{25jun121_XXX})-(\ref{25jun123_XXX}) for the case of homogeneous beliefs. The structure of the two systems is exactly the same, except that (\ref{25jun116})-(\ref{25jun123}) tracks a joint distribution for $(x,\theta)$ rather than a marginal distribution for $x$.

\subsection{Adaptive learning with common noise: sidestepping the Master equation}\label{sec:learning_no_noise_noise}
Starting from the second formulation in the preceding section with $(p,\theta)$ as state variables, the generalization to the case of common noise is straightforward. The key point of this section will be that, modeling adaptive learning in MFGs with a low-dimensional coupling and common noise, allows for sidestepping the infinite-dimensional Master equation.

Agents' perceived law of motion for prices is
\be\label{eq:general_learning_common}
d \widehat p_{s,t} = \mu_p(\widehat p_{s,t},Z_s,\theta)ds + \sigma_p(\widehat p_{s,t},Z_s,\theta)dB_s
\ee
where $\theta \in \Rm^d$ is a parameter vector. The agents' estimate of $\theta$ are still updated according to the learning rule \eqref{eq:general_learning_theta}.

Denoting the generator corresponding to \eqref{eq:general_learning_common} by $\mathcal{A}_p(z,\theta)$, we have the following HJB equation for the perceived value function $\widehat U(x,z,p,t;\theta)$ which is similar to \eqref{master_learning_no_noise}:
\be\label{eq:master_learning_noise}
\bal
\rho \widehat{U} - \partial_s \widehat{U} &= H(x,z,p,\nabla_x \widehat{U} ) 
+ \nu \Delta_x \widehat{U} + \mathcal{A}_p(z,\theta)\widehat U  + \beta \Delta_z \widehat{U}, \\
H(x,z,p,\lambda) &= \max_{\alpha \in A} \left\{ R(x,z,\alpha,p) + 
\alpha\cdot \lambda \right\},
\enbal
\ee
with  $(x,z,p,t)\in \ \Rm^n \times \Rm^k \times\Rm^\ell \times  (0,T), \theta \in \Rm^d$, and with corresponding policy
\be
\widehat \pi(x,z,p,t;\theta) = \arg \max_{\alpha \in A} \left\{ R(x,z,\alpha,p) + 
\alpha\cdot \nabla_x \widehat U(x,z,p,t;\theta) \right\}.
\ee
The important observation is that \eqref{eq:master_learning_noise} is a standard finite-dimensional HJB equation rather than an infinite-dimensional Master equation.

To solve for the evolution of the density $m_t$, equilibrium prices $p_t$, and the learned parameter estimates $\theta_t$, one can then proceed in the same fashion as in the case without common noise. First, define
\be\label{eq:policy_learning_noise}
\pi(x,Z_t,t)=\widehat\pi(x,Z_t,p_t,t;\theta_t)
\ee
where the price $p_t$ is given by
\be
p_t = P^*(m_t,Z_t).
\ee
Then solve the following forward-in-time system:
\be
\bal\label{system_learning_noise}
d m_t &= \mathcal{A}^*_{\pi_t,Z_t} m_t dt,\\
d Z_t &= \sqrt{2 \beta}dB_t\\
\dot{\theta}_t &= L(p_t, \theta_t)
\enbal
\ee
with the policy $\pi$ given by (\ref{eq:policy_learning_noise}) and with initial condition $(m_0,Z_0,\theta_0)$.

\paragraph{Remark on computational cost of \eqref{eq:master_learning_noise}.} The same remark as in the case without common noise applies: one does not actually have to solve \eqref{eq:master_learning_noise} for all values of $\theta \in \Rm^d$. This is because the equations for different $\theta$ are decoupled from each other. Starting from $\theta_0$, and simulating \eqref{system_learning_noise} forward in time, it is therefore sufficient to solve \eqref{eq:master_learning_noise} only for the $\theta$'s one actually encounters.

Either way, this computational cost should be compared to the cost of computing the infinite-dimensional Master equation with rational expectations. The former is clearly lower regardless of the exact computational strategy for \eqref{eq:master_learning_noise}.

\subsection{Other Directions: Reinforcement Learning and other Stochastic Approximation Methods}
One other approach for sidestepping the Master equation is to approximate the value function in Section \ref{sec:trouble} using ideas from the literature on reinforcement learning (RL). RL means learning value or policy functions of incompletely-known Markov decision processes via Monte Carlo simulation \citep{sutton-barto}. RL is typically formulated in discrete time but there exist continuous-time formulations \citep{doya,wang-zariphopoulou-zhou,jia-zhou}.

\citet{yang-etal-SRL} and \citet{wibault-etal} apply RL ideas to discrete-time MFGs with a low-dimensional coupling and common noise of the type studied here. Specifically, \citet{yang-etal-SRL} propose a ``Structural Policy Gradient" (SPG) method which leverages the known structure of individual dynamics (the generator $\mathcal{A}_\alpha$) while treating prices $p_t$ via simulation. While they impose that policy functions depend only on current prices,  \citet{wibault-etal} keep track of full price histories using recurrent neural networks as in recurrent RL \citep{hausknecht-stone,ni-etal-RNN}.

More generally, a promising approach could be to approximate these value functions using a ``stochastic approximation algorithm" \citep[e.g.][]{robbins-monro,ljung} of which reinforcement learning is a special case \citep{jaakkola-jordan-singh,tsitsiklis}.

\section{Discrete Time}\label{sec:discrete_time}
There is also a literature on discrete-time MFGs \citep[e.g.][]{gomes-mohr-riagosouza,gomes-mohr-riagosouza-2013} and this formulation may be useful for the application of some promising approaches to the challenge posed in this paper such as reinforcement learning. We therefore briefly repeat our paper's arguments in discrete time. For brevity we skip the case without common noise and focus directly on the more challenging case with common noise. For simplicity, we also focus on the case with finite state and action spaces, though this simplification is not essential. 

\subsection{Setup: Discrete-Time MFGs with Common Noise}
Everything is analogous to the continuous-time setup in Section \ref{sec:MFG}. Consider a system of $N\gg 1$ individual agents (players) at positions (states) $X_{i,t} \in \mathcal{X} \subset \Rm^n$, $i=1,\dots,N$, where $\mathcal{X}$ is a finite state space, i.e. $X_{i,t}$ can take only finitely many possible values. Time is discrete, $t=0,1,..,T$, where $T$ is a fixed terminal time that is sometimes taken as $T=+\infty$. As above, we consider the limit $N\to+\infty$ of a large number of agents. Anticipating this limit we here write equations directly in terms of the limiting density which we denote by $m_t(x) \in \mathcal{P}(\mathcal{X})$, the space of probability measures with support in $\mathcal{X}$. Given $x$ can take only finitely many values, this density is simply a high- but finite-dimensional vector (essentially a ``histogram"). The setup and notation are close to \citet{lauriere-etal-RL} and \citet{wibault-etal}.

Agents receive a period reward $R(X_{i,t},Z_t,\alpha_{i,t},m_t)$ and their state evolves according to a Markov process 
\begin{equation}\label{eq:x_markov}
X_{i,t+1} \sim \mathcal{T}_x(\cdot|X_{i,t},Z_t,\alpha_{i,t},m_t).
\end{equation}
As above, $\alpha_{i,t} \in A \subset \mathbb{R}^n$ is a control (but with $A$ a finite action space) and $Z_t \in \mathcal{Z} \subset \mathbb{R}^k$ is the aggregate state (common noise) that affects all agents and which evolves according to an exogenous Markov process:
\begin{equation}\label{eq:z_markov}
Z_{t+1} \sim \mathcal{T}_z(\cdot|Z_t).
\end{equation}
Agents maximize the cumulative discounted reward:
\be\label{eq:objective_discrete}
u_{i,0} = \max_{\alpha_i \in A} \mathbb{E} \left[\sum_{t=0}^T \gamma^t R(X_{i,t},Z_t,\alpha_{i,t},m_t) + \gamma^T V(X_{i,T},Z_T,m_T)\right]
\ee
subject to \eqref{eq:x_markov} and \eqref{eq:z_markov} and where $V$ is a terminal value and $0<\gamma \leq 1$ a discount factor.

As above, the optimal policy induces the density $m_t$ to evolve over time. This evolution is easiest to spell out for a slight generalization of the problem \eqref{eq:objective_discrete} in which we allow for stochastic policies of the form:
\be\label{eq:policy_stochastic}
\alpha_{i,t} \sim \pi_t(\cdot|X_{i,t},Z_t,m_t)
\ee
where $\pi_t$ is the probability distribution over actions $\alpha_{i,t}$ conditional on the states $(X_{i,t},Z_t,m_t)$. Given the optimal policy $\pi_t$, the density $m_t$ then evolves according to a Chapman-Kolmogorov equation (the discrete-time analogue of a Fokker-Planck equation)
\begin{equation}\label{eq:chapman}
m_{t+1}(x) = \sum_{\tilde{x},\tilde{\alpha}} m_t(\tilde{x}) \pi_t(\tilde{\alpha}|\tilde{x},Z_t,m_t)\mathcal{T}_x(x|\tilde{x},\tilde{\alpha},Z_t,m_t)
\end{equation}
which we can also write
\be \label{eq:chapman_matrix}
m_{t+1} = \mathbf{A}_{\pi_t,Z_t}^{\rm T}m_t,
\ee
where $\mathbf{A}_{\pi_t,Z_t}$ is the transition matrix of $x$ induced by the optimal policy $\pi_t$. This equation is the discrete-time counterpart to the Fokker-Planck equation \eqref{eq:dynamics_RE}, with the transition matrix $\mathbf{A}_{\pi_t,Z_t}$ being the counterpart to the generator $\mathcal{A}_{\pi,Z_t}$ (and the matrix transpose that to the operator adjoint). Note that, because $Z_t$ evolves according to the Markov process \eqref{eq:z_markov}, the density $m_{t+1}$ is itself a stochastic process.

\paragraph{Mean Field Games with a low-dimensional coupling.} As in Section \ref{sec:special_structure}, we pay special attention to the class of ``MFGs with a low-dimensional coupling." Using analogous notation, agents choose policies $\pi$ to maximize
\begin{equation}\label{objective_MFG_lowD_discrete}
u_{i,0} = \max_{\alpha_i \in A}\mathbb{E} \left[\sum_{t=0}^T \gamma^t \widetilde R(X_{i,t},Z_t,\alpha_{i,t},p_t) + \gamma^T \widetilde V(X_{i,T},Z_T,p_T)\right]
\end{equation}
subject to 
\be \label{eq:policy_discrete_lowD}
X_{i,t+1} \sim \mathcal{T}_x(\cdot|X_{i,t},Z_t,\alpha_{i,t},p_t), 
\ee
and \eqref{eq:z_markov} where
\be\label{eq:functional_discrete}
p_t= P^*(m_t,Z_t),
\ee
for a fixed functional $P^*: \mathcal{P}(\mathcal{X}) \times \mathcal{Z} \rightarrow \mathbb{R}^\ell$. Again note that model agents do not directly ``care about" the density $m_t$ in the sense that it does not enter their running reward functions; instead they only ``care about" the much lower-dimensional vector $p_t$. As discussed in Sections \ref{sec:macroMFG} and \ref{sec:macro_example} such MFGs arise naturally in macroeconomics where they are known as ``heterogeneous agent models" \citep[e.g.][]{krusell-smith,denhaan}. In this case, the vector $p_t$ has the interpretation of ``equilibrium prices" that are determined from some ``market clearing" conditions analogous to \eqref{25apr2702}.

\paragraph{Rational expectations.}  As above, standard formulations of this MFG impose rational expectations: agents know not only the correct transition probabilities $\mathcal{T}_x$ and $\mathcal{T}_z$ but also the high-dimensional transition matrix $\mathbf{A}_{\pi_t,Z_t}$ which governs the evolution of the complex system they inhabit.

\paragraph{Master equation in discrete-time MFGs.}
In both general MFGs and MFGs with a low-dimensional coupling, under rational expectations, the agents' optimization problem gives rise to the Master equation (i.e. Bellman equation on the space of probability measures):
\begin{equation}\label{eq:master_discrete}
\begin{split}
U_t(x,z,m) &= \max_\alpha \ R(x,z,\alpha,m) + \gamma \mathbb{E}_{x',z'}[U_{t+1}(x',z',m')|x,z,m] \quad \mbox{subject to}\\
x' &\sim \mathcal{T}_x(\cdot|x,z,\alpha,m),\\
z' &\sim \mathcal{T}_z(\cdot|z),\\
m'&= \mathbf{A}_{\pi_t,z}^{\rm T} m,\\
U_T(x,z,m) &= V(x,z,m).
\end{split}
\end{equation}
This equation is the exact discrete-time counterpart to \eqref{eq:master_RE} and it similarly features the state variable $m \in \mathcal{P}(\mathcal{X})$, the space of probability measures with support in $\mathcal{X}$. It therefore suffers from an extreme version of the curse of dimensionality.

Section \ref{sec:trouble} criticized the continuous-time Master equation in MFGs with a low-dimensional coupling. All the same criticisms also apply to its discrete-time counterpart.

%

\subsection{Discrete-Time MFGs without Rational Expectations}\label{sec:nonRE_discrete}

\paragraph{General MFGs without rational expectations.} As above, we assume that agents have rational expectations about the evolution of their own individual state $X_{i,t}$, i.e. that they know the correct transition probabilities $\mathcal{T}_x$. However, like in Section \ref{sec:nonRE_MFG_common}, we allow for the possibility that they may have non-rational expectations about the evolution of the density $m_t$ and the aggregate state $Z_t$. In particular, agents believe that the density and aggregate state evolve according to the perceived laws of motion
\be \label{eq:chapman_matrix_nonRE}
\bal 
\widehat{m}_{s+1,t} &= \mathbf{B}_{\widehat Z_t}^{\rm T} \widehat{m}_{s,t}, \quad s \geq t\\
\widehat Z_{s+1,t} &\sim \widehat{\mathcal{T}}_z(\cdot|\widehat Z_{s,t}), \quad s \geq t,
\enbal
\ee
with $\widehat{m}_{t,t}=m_t$ and $\widehat Z_{t,t}=Z_t$ rather than \eqref{eq:z_markov} and \eqref{eq:chapman_matrix}. Under this assumption, agents' optimization problem gives rise to a Master equation that is just like \eqref{eq:master_discrete} but with $\mathcal{T}_z$ replaced by $\widehat \mathcal{T}_z$ and $\mathbf{A}_{\pi,z}$ replaced by $\mathbf{B}_{\widehat z}$. Rational expectations imposes that
$$\mathbf{B}_{z} = \mathbf{A}_{\pi,z} \quad \mbox{and} \quad \widehat{\mathcal{T}}_z = \mathcal{T}_z $$
so that, in the special case of rational expectations, we recover the usual Master equation \eqref{eq:master_discrete}.

\paragraph{Sidestepping the Master equation in MFGs with a low-dimensional coupling.} It is again most interesting to consider non-rational expectations in MFGs with a low-dimensional coupling because this holds the promise of sidestepping the Master equation altogether.

With rational expectations, agents understand the dependence of $p_t$ on $m_t$ and $Z_t$ and therefore use the function $P^*$ together with the correct stochastic processes for $m_t$ and $Z_t$ to predict future values of $p_t$. This leads to the Master equation (\ref{eq:master_discrete}), with essentially no simplifications despite the low-dimensional coupling.

With non-rational expectations, agents instead perceive some other stochastic process for the pair~$(p_t,Z_t)$. Like in Section  \ref{sec:sidestepping_lowD}, they could simply perceive $p_t$ to evolve according to an exogenous Markov process
\begin{equation}\label{eq:price_process_simple_discrete}
\widehat{p}_{s+1,t} \sim \widehat \mathcal{T}_p(\cdot|\widehat{p}_{s,t}), \quad s \geq t, \quad \widehat{p}_{t,t} = p_t.
\end{equation}
In more complicated cases, agents may perceive a joint stochastic process for $p_t,Z_t$ and other variables.

In the case of agents perceiving the simple process \eqref{eq:price_process_simple_discrete}, instead of writing a Master equation, we can write a much simpler, 
standard finite-dimensional Bellman equation (as above, going forward, we drop the tildes from $\widetilde R$ and $\widetilde V$ for notational simplicity):
\begin{equation}\label{eq:master_nonRE_CE_discrete}
\begin{split}
\widehat U_t(x,z,p) &= \max_\alpha \ R(x,z,\alpha,p) + \gamma \mathbb{E}_{x',z',p'}[U_{t+1}(x',z',p')|x,z,p] \quad \mbox{subject to}\\
x' &\sim \mathcal{T}_x(\cdot|x,z,\alpha,p),\\
z' &\sim \mathcal{T}_z(\cdot|z),\\
p' &\sim \widehat \mathcal{T}_p(\cdot|p), \\
\widehat U_T(x,z,p) &= V(x,z,p).
\end{split}
\end{equation}
Therefore, in MFGs with a low-dimensional coupling, departing from rational expectations can completely sidestep the Master equation. Of course, the case considered here is just an illustrative example. In particular, note that the perceived law of motion \eqref{eq:price_process_simple} is specified completely ``outside the model" which leaves open the question where this perceived law of motion ``comes from" in the first place.

\subsection{The Challenge and a Markov Reward Process}
As discussed in Section \ref{sec:learning} and \citet{moll-challenge}, the challenge is: how can we formulate, in a systematic way, models of agents' behavior in situations with a low-dimensional coupling that lead to equations that (i)~approximate agents' real-world behavior, and (ii) sidestep computing the solutions to a Master equation with the infinite-dimensional state $m \in \mathcal{P}(\mathcal{X})$ and the associated curse of dimensionality?

\paragraph{A Markov Reward Process with all the difficulty.}
To understand the key difficulty, it is useful to consider a simplified version of the model with no actions $\alpha$, i.e. a Markov Reward Process (MRP) rather than a Markov Decision Process (MDP). This material is adapted from ongoing work by \citet{yang-etal-SRL}.

To this end, eliminate actions $\alpha$ and assume that the state $X_{i,t}$ evolves exogenously according to a transition matrix $\mathbf{A}_{Z_t}$ that depends on the aggregate state $Z_t$, implying that $(m_t,Z_t)$ evolve as
\begin{equation}\label{eq:transition_MRP}
m_{t+1} = \mathbf{A}_{Z_t}^{\rm T}m_t, \quad Z_{t+1} \sim \mathcal{T}_z(\cdot|Z_t).
\end{equation}
Similarly, replace the running reward and terminal value in \eqref{objective_MFG_lowD_discrete}  by a reward $R(p)$ and terminal value $V(p)$ that depend only on the low-dimensional price vector $p_t \in \mathbb{R}^\ell$. As before, the vector $p_t$ is still determined by \eqref{eq:functional_discrete}.

The simplified problem is to compute the value of a MRP: given $p_t$ evolving according to \eqref{eq:functional_discrete} and $(m_t,Z_t)$ according to \eqref{eq:transition_MRP}, compute the expected presented discounted value (PDV) of rewards:
\begin{equation}\label{objective_MRP}
u_{0} = \mathbb{E} \left[\sum_{t=0}^T \gamma^t R(p_t) + \gamma^T V(p_T)\right].
\end{equation}
This problem contains all the difficulty of the more complicated problem in MFGs with a low-dimensional coupling.

The ``correct" way -- in the rational expectations sense -- of computing the value of this MRP is to solve a Master equation for the value function $U_t(z,m)$:
\begin{equation}\label{eq:master_MRP}
\begin{split}
U_t(z,m) &= R(P^*(m,z)) + \gamma \mathbb{E}_{z'}[U_{t+1}(z',m')|z,m] \quad \mbox{subject to}\\
z' &\sim \mathcal{T}_z(\cdot|z),\\
m'&= \mathbf{A}_{z}^{\rm T} m\\
U_T(z,m) &= V(P^*(m,z))
\end{split}
\end{equation}
This Master equation illustrates the trouble with the rational expectations assumption: even though the reward function is only a function of a low-dimensional vector  $p_t \in \mathbb{R}^\ell$, computing the PDV in \eqref{objective_MRP} requires solving a Bellman equation on the space of probability measures $\mathcal{P}(\mathcal{X})$.

\paragraph{The difficulty: prices $p_t$ are not Markov.} As discussed in Section \ref{sec:nonRE_lowD}, the key difficulty is that the vector $p_t$ does not follow a Markov process; instead only $(m_t,Z_t)$ has the Markov property and $p_t$ is instead a complicated non-linear functional of this Markov state. Agents with rational expectations therefore (unrealistically) forecast the Markov state $(m_t,Z_t)$ in order to forecast the non-Markovian $p_t$.

\subsection{Adaptive Learning in Discrete Time}
Finally, we show how to write the adaptive learning model of Section \ref{sec:learning} in discrete time. The economics literature typically formulates such models in discrete time \citep{bray,marcet-sargent,evans-honkapohja,jacobson}. The key assumption is that, to forecast prices, agents use a \emph{perceived law of motion} in the form of a Markov process:
\begin{equation}\label{eq:price_process_learning_discrete}
\widehat{p}_{s+1,t} \sim \widehat \mathcal{T}_p(\cdot|\widehat{p}_{s,t},Z_s,\theta), \quad s \geq t, \quad \widehat{p}_{t,t} = p_t,
\end{equation}
where $\theta \in \Rm^d$ is a parameter vector. The key difference to the Markov process \eqref{eq:price_process_simple_discrete} in Section \ref{sec:nonRE_discrete} is that the process is endogenous to the model because agents learn the parameter vector $\theta$ over time from past observations of $p_t$. Specifically, they form an estimate $\widehat{\theta}_{t}$ of $\theta$ which they update using the learning rule
\be\label{eq:learning_discrete}
\widehat{\theta}_{t+1} = L(p_t,\widehat{\theta}_t).
\ee
For example, \eqref{eq:price_process_learning_discrete} could be a vector autoregressive (VAR) process for the vector of prices and \eqref{eq:learning_discrete} a recursive least-squares estimator for the parameters of this VAR.

Dropping the hat subscripts from $\widehat{\theta}_t$ for notational simplicity (but keeping in 
mind that this is really a time-varying estimate of the parameter $\theta$), agents' optimization problem given a current estimate $\theta$ is:
\be
\widehat{U}_t(x,z,p;\theta) = \max_{\alpha_i \in A} \mathbb{E} 
\left[\sum_{s=t}^T \gamma^{s-t} R(X_{i,s},Z_t,\alpha_{i,s},\widehat p_{s,t}) + \gamma^{T-t} V(X_{i,T},Z_T,\widehat p_{T,t}) \right]
\ee
subject to
\begin{align}
\widehat{p}_{s+1,t} &\sim \widehat \mathcal{T}_p(\cdot|\widehat{p}_{s,t},Z_s,\theta), \quad s \geq t,\\
\widehat{p}_{t,t} &= p
\end{align}
and subject to \eqref{eq:z_markov} and \eqref{eq:policy_discrete_lowD}.

The corresponding Bellman equation is:
\begin{equation}\label{eq:master_learning_discrete}
\begin{split}
\widehat U_t(x,z,p;\theta) &= \max_\alpha \ R(x,z,\alpha,p) + \gamma \mathbb{E}_{x',z',p'}[U_{t+1}(x',z',p';\theta)|x,z,p] \quad \mbox{subject to}\\
x' &\sim \mathcal{T}_x(\cdot|x,z,\alpha,p),\\
z' &\sim \mathcal{T}_z(\cdot|z),\\
p' &\sim \widehat \mathcal{T}_p(\cdot|p,z,\theta), \\
\widehat U_T(x,z,p;\theta) &= V(x,z,p),
\end{split}
\end{equation}
with corresponding optimal policy $\widehat \pi_t(x,z,p,\theta)$. The key observation is that this is a standard finite-dimensional Bellman equation rather than an infinite-dimensional Master equation. 
Analogous to the discussion in Section \ref{sec:learning_no_noise}, the Bellman equation \eqref{eq:master_learning_discrete} assumes that learning is external to agents and one can instead formulate a variant with internalized learning \citep[e.g.][]{christiano-eichenbaum-johannsen}. To obtain this variant, we simply replace $\theta$ in $U_{t+1}(x',z',p';\theta)$ by $\theta'$ and add the law of motion $\theta' =L(p,\theta)$ to the constraint set.

With the solution in hand, the evolution of the density $m_t$ and equilibrium prices $p_t$ are found as follows. First, define
\be\label{eq:policy_learning_noise_discrete}
\pi_t(x,Z_t)=\widehat\pi_t(x,Z_t,p_t;\theta_t)
\ee
where the price $p_t$ is given by
\be
p_t = P^*(m_t,Z_t).
\ee
Then solve the following forward-in-time system:
\be
\bal\label{system_learning_noise_discrete}
m_{t+1} &= \mathbf{A}_{\pi_t,Z_t}^{\rm T}m_t\\
Z_{t+1} &\sim \mathcal{T}_z(\cdot|Z_t)\\
\theta_{t+1} &= L(p_t, \theta_t)
\enbal
\ee
with the policy $\pi_t$ given by (\ref{eq:policy_learning_noise_discrete}) and with initial condition $(m_0,Z_0,\theta_0)$. A similar problem is solved by \citet{jacobson}.

\paragraph{Remark on computational cost of \eqref{eq:master_learning_discrete}.} As noted in Section \ref{sec:learning}, one does not actually have to solve \eqref{eq:master_learning_discrete} for all values of $\theta \in \Rm^d$. This is because the equations for different $\theta$ are decoupled from each other. Starting from $\theta_0$, and simulating \eqref{system_learning_noise_discrete} forward in time, it is therefore sufficient to solve \eqref{eq:master_learning_discrete} only for the $\theta$'s one actually encounters.

\paragraph{Remark on relation to \citet{krusell-smith} algorithm.} There is a link between this adaptive learning approach (specifically, the variant with least-squares learning) and the algorithm of \citet{krusell-smith}. In both approaches, decision makers use a perceived law of motion like \eqref{eq:price_process_learning_discrete} and estimate its coefficients via least squares. A difference is that adaptive learning updates the coefficient estimate $\theta_{t+1}$ from $\theta_t$ recursively over time so that solving for the MFG equilibrium and belief updating are done ``in one sweep" via solving \eqref{system_learning_noise_discrete} forward in time.

\section{Conclusion}\label{sec:conclusion}
This paper has shown how to formulate MFGs without rational expectations, i.e. without the assumption that agents know all relevant transition probabilities for the complex system they inhabit. Instead of using the correct transition probabilities, agents instead use some other ``non-rational" transition probabilities when solving their optimization problems. We show how to write the corresponding equations describing the Nash equilibrium of the MFG, both for the case with and without common noise. In the special case of rational expectations we recover the standard backward-forward MFG system and MFG Master equation.

Departing from rational expectations is particularly relevant when there is common noise in ``MFGs with a low-dimensional coupling", i.e. MFGs in which agents' running reward function depends on the density only through low-dimensional functionals and which are typical in macroeconomics. In MFGs with a low-dimensional coupling, departing from rational expectations allows for completely sidestepping the Master equation and for instead solving finite-dimensional HJB equations. We introduced an adaptive learning model as a particular example of non-rational expectations and discussed its properties.

\small
\setstretch{1.0}{\bibliographystyle{aer}
\bibliography{challenge_bib}

@Article{marcet-sargent,
  author={Marcet, Albert and Sargent, Thomas J.},
  title={{Convergence of Least Squares Learning Mechanisms in Self-referential Linear Stochastic Models}},
  journal={Journal of Economic Theory},
  year=1989,
  volume={48},
  number={2},
  pages={337-368},
  month={August},
  doi={},
  url={https://ideas.repec.org/a/eee/jetheo/v48y1989i2p337-368.html}
}

@Article{aiyagari,
  author={Aiyagari, S. Rao},
  title={{Uninsured Idiosyncratic Risk and Aggregate Saving}},
  journal={The Quarterly Journal of Economics},
  year=1994,
  volume={109},
  number={3},
  pages={659-84},
  month={August}
}

@Article{krusell-smith,
  author={Per Krusell and Anthony A. Smith},
  title={{Income and Wealth Heterogeneity in the Macroeconomy}},
  journal={Journal of Political Economy},
  year=1998,
  volume={106},
  number={5},
  pages={867-896},
  month={October}
}

@Article{denhaan,
  author={{Den Haan}, Wouter J.},
  title={{Heterogeneity, Aggregate Uncertainty, and the Short-Term Interest Rate}},
  journal={Journal of Business \& Economic Statistics},
  year=1996,
  volume={14},
  number={4},
  pages={399-411},
  month={October}
}

@article{AHLLM,
    author = {Achdou, Yves and Han, Jiequn and Lasry, Jean-Michel and Lions, Pierre-Louis and Moll, Benjamin},
    title = "{{Income and Wealth Distribution in Macroeconomics: A Continuous-Time Approach}}",
    journal = {The Review of Economic Studies},
    volume = {89},
    number = {1},
    pages = {45-86},
    year = {2021},
    month = {04},
    issn = {0034-6527},
    doi = {10.1093/restud/rdab002},
    url = {https://doi.org/10.1093/restud/rdab002},
    eprint = {https://academic.oup.com/restud/article-pdf/89/1/45/42137446/rdab002.pdf},
}

@Article{grandmont,
  author={Grandmont, Jean-Michel},
  title={{Temporary General Equilibrium Theory}},
  journal={Econometrica},
  year=1977,
  volume={45},
  number={3},
  pages={535-572},
  month={April}
}

@book{hicks-value,
  author={John R. Hicks},
  title={{Value and Capital: An Inquiry into Some Fundamental Principles of Economic Theory}},
  year=1939,
  publisher = {Clarendon Press. Available at \url{https://benjaminmoll.com/Hicks_Value_and_Capital/}},
  url={https://benjaminmoll.com/Hicks_Value_and_Capital/}
}

@book{lindahl,
  author={Erik Lindahl},
  title={{Studies in the Theory of Money and Capital}},
  year=1939,
  publisher = {Allen and Unwin}
}

@TechReport{moll-challenge,
  author={Benjamin Moll},
  title={{The Trouble with Rational Expectations in Heterogeneous Agent Models: A Challenge for Macroeconomics}},
  year=2025,
  institution={{Royal Economic Society, \url{https://benjaminmoll.com/challenge/}}},
  type={{Economic Journal Lecture}},
  url={https://benjaminmoll.com/challenge/}
}

@InProceedings{xu-etal-RL,
  title =    {{Finding Regularized Competitive Equilibria of Heterogeneous Agent Macroeconomic Models via Reinforcement Learning}},
  author =       {Xu, Ruitu and Min, Yifei and Wang, Tianhao and Jordan, Michael I. and Wang, Zhaoran and Yang, Zhuoran},
  booktitle =    {Proceedings of The 26th International Conference on Artificial Intelligence and Statistics},
  pages =    {375--407},
  year =   {2023},
  editor =   {Ruiz, Francisco and Dy, Jennifer and van de Meent, Jan-Willem},
  volume =   {206},
  series =   {Proceedings of Machine Learning Research},
  month =    {25--27 Apr},
  publisher =    {PMLR},
  url =    {https://proceedings.mlr.press/v206/xu23a.html}
}

@book{sutton-barto,
  title={{Reinforcement Learning: An Introduction}},
  author={Sutton, Richard S. and Barto, Andrew G.},
  year={2018},
  publisher={MIT Press}
}

@article{muth,
  title={{Rational Expectations and the Theory of Price Movements}},
  author={Muth, John F.},
  journal={Econometrica},
  pages={315--335},
  year={1961},
 volume = {29},
number = {3},
  month={July}
}

@article{robbins-monro,
 ISSN = {00034851},
 URL = {http://www.jstor.org/stable/2236626},
 author = {Herbert Robbins and Sutton Monro},
 journal = {The Annals of Mathematical Statistics},
 number = {3},
 pages = {400--407},
 publisher = {Institute of Mathematical Statistics},
 title = {{A Stochastic Approximation Method}},
 urldate = {2024-07-19},
 volume = {22},
 year = {1951}
}

@article{tsitsiklis,
  title={{Asynchronous Stochastic Approximation and Q-learning}},
  author={Tsitsiklis, John N.},
  journal={Machine Learning},
  volume={16},
  pages={185--202},
  year={1994}
}

@article{jaakkola-jordan-singh,
  title={{Convergence of Stochastic Iterative Dynamic Programming Algorithms}},
  author={Jaakkola, Tommi and Jordan, Michael and Singh, Satinder},
  journal={Advances in Neural Information Processing Systems},
  volume={6},
  year={1993},
  pages={703--710}

}

@InProceedings{lauriere-etal-RL,
  title = 	 {{Scalable Deep Reinforcement Learning Algorithms for Mean Field Games}},
  author =       {Lauri\`{e}re, Mathieu and Perrin, Sarah and Girgin, Sertan and Muller, Paul and Jain, Ayush and Cabannes, Theophile and Piliouras, Georgios and Perolat, Julien and Elie, Romuald and Pietquin, Olivier and Geist, Matthieu},
  booktitle = 	 {Proceedings of the 39th International Conference on Machine Learning},
  pages = 	 {12078--12095},
  year = 	 {2022},
  editor = 	 {Chaudhuri, Kamalika and Jegelka, Stefanie and Song, Le and Szepesvari, Csaba and Niu, Gang and Sabato, Sivan},
  volume = 	 {162},
  series = 	 {Proceedings of Machine Learning Research},
  month = 	 {17--23 Jul},
  publisher =    {PMLR},
  url = 	 {https://proceedings.mlr.press/v162/lauriere22a.html}
}

@misc{lauriere-etal-survey,
      title={Learning in Mean Field Games: A Survey}, 
      author={Mathieu Lauri\`{e}re and Sarah Perrin and Julien Pérolat and Sertan Girgin and Paul Muller and Romuald Élie and Matthieu Geist and Olivier Pietquin},
      year={2024},
      eprint={2205.12944},
      archivePrefix={arXiv},
      primaryClass={cs.LG},
      url={https://arxiv.org/abs/2205.12944}, 
}

@Article{bray,
  author={Bray, Margaret},
  title={{Learning, Estimation, and the Stability of Rational Expectations}},
  journal={Journal of Economic Theory},
  year=1982,
  volume={26},
  number={2},
  pages={318-339},
  month={April},
  doi={},
  url={https://ideas.repec.org/a/eee/jetheo/v26y1982i2p318-339.html}
}

@ARTICLE{ljung,
  author={Ljung, L.},
  journal={IEEE Transactions on Automatic Control}, 
  title={{Analysis of Recursive Stochastic Algorithms}}, 
  year={1977},
  volume={22},
  number={4},
  pages={551-575},
  doi={10.1109/TAC.1977.1101561}
}

@book {cardaliaguet-delarue-lasry-lions,
    AUTHOR = {Cardaliaguet, Pierre and Delarue, Fran\c{c}ois and Lasry,
              Jean-Michel and Lions, Pierre-Louis},
     TITLE = {{The Master Equation and the Convergence Problem in Mean Field
              Games}},
    SERIES = {Annals of Mathematics Studies},
    VOLUME = {201},
 PUBLISHER = {Princeton University Press},
      YEAR = {2019},
     PAGES = {x+212},
      ISBN = {978-0-691-19071-6; 978-0-691-19070-9},
   MRCLASS = {49-02 (49N70 60H30 60K35 91A13 91A15)},
  MRNUMBER = {3967062},
       DOI = {10.2307/j.ctvckq7qf},
       URL = {https://doi.org/10.2307/j.ctvckq7qf},
}

@article {lasry-lions,
   author = {Lasry, Jean-Michel and Lions, Pierre-Louis},
   title = {Mean field games},
   journal = {Japanese Journal of Mathematics},
   publisher = {Springer Japan},
   issn = {0289-2316},
   pages = {229-260},
   volume = {2},
   issue = {1},
   year = {2007}
}

@incollection{grandmont-palgrave,
  title={{Temporary Equilibrium}},
  author={Grandmont, Jean-Michel},
  booktitle={General Equilibrium},
  pages={297--304},
  year={1989},
  publisher={Springer}
}

@techreport{christiano-eichenbaum-johannsen,
 title = {{Slow Learning}},
 institution = "National Bureau of Economic Research",
author = {Christiano, Lawrence J. and Eichenbaum, Martin S. and Johannsen, Benjamin K.},
 type = "NBER Working Papers",
 number = "32358",
 year = "2024",
 month = "April",
 doi = {10.3386/w32358},
 URL = "http://www.nber.org/papers/w32358",
 
}

@book{evans-honkapohja,
  title={{Learning and Expectations in Macroeconomics}},
  author={George W. Evans and Seppo Honkapohja},
  year={2001},
  publisher={Princeton University Press},
}

@article{lucas-critique,
title = {{Econometric Policy Evaluation: A Critique}},
journal = {Carnegie-Rochester Conference Series on Public Policy},
volume = {1},
pages = {19-46},
year = {1976},
issn = {0167-2231},
doi = {https://doi.org/10.1016/S0167-2231(76)80003-6},
url = {https://www.sciencedirect.com/science/article/pii/S0167223176800036},
author = {Robert E. Lucas}
}

@misc{cardaliaguet,
  author={Pierre Cardaliaguet},
  title={{Notes on Mean-Field Games (from P.-L. Lions' Lectures at Coll\`{e}ge de France)}},
  howpublished = {\url{https://www.ceremade.dauphine.fr/~cardaliaguet/MFG20130420.pdf}},
  year=2013,
  month= {},
  institution={Dauphine},
}

@misc{ryzhik,
  author={Lenya Ryzhik},
  title={{Lecture Notes for a Reading Course on Mean-Field Games}},
  howpublished = {\url{https://math.stanford.edu/~ryzhik/STANFORD/MEAN-FIELD-GAMES/notes-mean-field.pdf}},
  year=2018,
  month= {},
  institution={Stanford},
}

@incollection{bewley,
    address = {Amsterdam},
    author = {Bewley, Truman},
    booktitle = {Contributions to Mathematical Economics in Honor of Gerard Debreu},
    editor = {Hildenbrand, Werner and Mas-Collel, Andreu},
    publisher = {North-Holland},
    title = {{Stationary Monetary Equilibrium with a Continuum of Independently Fluctuating Consumers}},
    year = {1986}
}

@Article{huggett,
  author={Huggett, Mark},
  title={The risk-free rate in heterogeneous-agent incomplete-insurance economies},
  journal={Journal of Economic Dynamics and Control},
  year=1993,
  volume={17},
  number={5-6},
  pages={953-969},
  month={}
}

@article{PDE-macro,
author = {Achdou, Yves  and Buera, Francisco J.  and Lasry, Jean-Michel  and Lions, Pierre-Louis  and Moll, Benjamin },
title = {Partial differential equation models in macroeconomics},
journal = {Philosophical Transactions of the Royal Society A: Mathematical, Physical and Engineering Sciences},
volume = {372},
number = {2028},
pages = {20130397},
year = {2014},
doi = {10.1098/rsta.2013.0397},
URL = {https://royalsocietypublishing.org/doi/abs/10.1098/rsta.2013.0397},
eprint = {https://royalsocietypublishing.org/doi/pdf/10.1098/rsta.2013.0397}
}

@article{doya,
    author = {Doya, Kenji},
    title = {Reinforcement Learning in Continuous Time and Space},
    journal = {Neural Computation},
    volume = {12},
    number = {1},
    pages = {219-245},
    year = {2000},
    month = {01}
}

@article{wang-zariphopoulou-zhou,
  author  = {Haoran Wang and Thaleia Zariphopoulou and Xun Yu Zhou},
  title   = {Reinforcement Learning in Continuous Time and Space: A Stochastic Control Approach},
  journal = {Journal of Machine Learning Research},
  year    = {2020},
  volume  = {21},
  number  = {198},
  pages   = {1--34},
  url     = {http://jmlr.org/papers/v21/19-144.html}
}

@book {carmona-delarue1,
    AUTHOR = {Carmona, Ren\'{e} and Delarue, Fran\c{c}ois},
     TITLE = {Probabilistic theory of mean field games with applications.
              {I}},
    SERIES = {Probability Theory and Stochastic Modelling},
    VOLUME = {83},
      NOTE = {Mean field FBSDEs, control, and games},
 PUBLISHER = {Springer, Cham},
      YEAR = {2018},
     PAGES = {xxv+713},
      ISBN = {978-3-319-56437-1; 978-3-319-58920-6},
   MRCLASS = {60-02 (35R60 49N70 49N90 60H15 60H30 91A15 93E20)},
  MRNUMBER = {3752669},
MRREVIEWER = {Vassili N. Kolokol\cprime tsov},
}

@book {carmona-delarue2,
    AUTHOR = {Carmona, Ren\'{e} and Delarue, Fran\c{c}ois},
     TITLE = {Probabilistic theory of mean field games with applications.
              {II}},
    SERIES = {Probability Theory and Stochastic Modelling},
    VOLUME = {84},
      NOTE = {Mean field games with common noise and master equations},
 PUBLISHER = {Springer, Cham},
      YEAR = {2018},
     PAGES = {xxiv+697},
      ISBN = {978-3-319-56435-7; 978-3-319-56436-4},
   MRCLASS = {60-02 (35R60 49L20 60G55 60H10 60H30 91A13 91A15)},
  MRNUMBER = {3753660},
MRREVIEWER = {Vassili N. Kolokol\cprime tsov},
}

@article{jia-zhou,
author = {Jia, Yanwei and Zhou, Xun Yu},
title = {Q-learning in continuous time},
year = {2023},
issue_date = {January 2023},
publisher = {JMLR.org},
volume = {24},
number = {1},
issn = {1532-4435},
journal = {J. Mach. Learn. Res.},
month = jan,
articleno = {161},
numpages = {61}
}

@TechReport{jacobson,
  author={Margaret M. Jacobson},
  title={{Beliefs, Aggregate Risk, and the U.S. Housing Boom}},
  year=2025,
  month=Sep,
  institution={Board of Governors of the Federal Reserve System (U.S.)},
  type={Finance and Economics Discussion Series},
  url={https://ideas.repec.org/p/fip/fedgfe/2022-61.html},
  number={2022-061}
}

@misc{yang-etal-SRL,
      title={Structural Reinforcement Learning for Heterogeneous Agent Macroeconomics}, 
      author={Yucheng Yang and Chiyuan Wang and Andreas Schaab and Benjamin Moll},
      year={2025},
      eprint={2512.18892},
      archivePrefix={arXiv},
      primaryClass={econ.TH},
      url={https://arxiv.org/abs/2512.18892}, 
}

@misc{wibault-etal,
      title={{Recurrent Structural Policy Gradient for Partially Observable Mean Field Games}}, 
      author={Clarisse Wibault and Sebastian Towers and Tiphaine Wibault and Juan Duque and Johannes Forkel and George Whittle and
Andreas Schaab and Yucheng Yang and Chiyuan Wang and Michael Osborne and Benjamin Moll and Jakob Foerster},
      year={2026}
}

@misc{hausknecht-stone,
      title={Deep Recurrent Q-Learning for Partially Observable MDPs}, 
      author={Matthew Hausknecht and Peter Stone},
      year={2017},
      eprint={1507.06527},
      archivePrefix={arXiv},
      primaryClass={cs.LG},
      url={https://arxiv.org/abs/1507.06527}, 
}

@misc{ni-etal-RNN,
      title={Recurrent Model-Free RL Can Be a Strong Baseline for Many POMDPs}, 
      author={Tianwei Ni and Benjamin Eysenbach and Ruslan Salakhutdinov},
      year={2022},
      eprint={2110.05038},
      archivePrefix={arXiv},
      primaryClass={cs.LG},
      url={https://arxiv.org/abs/2110.05038}, 
}

@article{simon-rationality,
 ISSN = {00028282},
 URL = {http://www.jstor.org/stable/1816653},
 author = {Herbert A. Simon},
 journal = {The American Economic Review},
 number = {2},
 pages = {1--16},
 publisher = {American Economic Association},
 title = {{Rationality as Process and as Product of Thought}},
 urldate = {2025-07-18},
 volume = {68},
 year = {1978}
}

@incollection{cardal-porr,
  author    = {Cardaliaguet, Pierre and Porretta, Alessio},
  title     = {An Introduction to Mean Field Game Theory},
  booktitle = {Mean Field Games},
  series    = {Lecture Notes in Mathematics},
  volume    = {2281},
  pages     = {1--158},
  publisher = {Springer},
  address   = {Cham},
  year      = {2020}
}

@unpublished{bertucci-meynard1,
  author = {Bertucci, Charles and Meynard, Charles},
  title  = {Noise through an Additional Variable for Mean Field Games Master Equation on Finite State Space},
  year   = {2024},
  note   = {Working paper}
}

@unpublished{bertucci-meynard2,
  author = {Bertucci, Charles and Meynard, Charles},
  title  = {A Study of Common Noise in Mean Field Games},
  year   = {2024},
  note   = {Working paper}
}

@article{ahuja-2016,
  author  = {Ahuja, Saran},
  title   = {Wellposedness of Mean Field Games with Common Noise under a Weak Monotonicity Condition},
  journal = {SIAM Journal on Control and Optimization},
  volume  = {54},
  number  = {1},
  pages   = {30--48},
  year    = {2016}
}

@article{ahuja-2019,
  author  = {Ahuja, Saran and Ren, Weiluo and Yang, Tzu-Wei},
  title   = {Forward-Backward Stochastic Differential Equations with Monotone Functionals and Mean Field Games with Common Noise},
  journal = {Stochastic Processes and their Applications},
  volume  = {129},
  number  = {10},
  pages   = {3859--3892},
  year    = {2019}
}

@unpublished{bertucci-incomplete,
  author = {Bertucci, Charles},
  title  = {Mean Field Games with Incomplete Information},
  year   = {2023},
  note   = {Working paper}
}

@incollection{cagan,
  author    = {Cagan, Phillip D.},
  title     = {The Monetary Dynamics of Hyperinflation},
  booktitle = {Studies in the Quantity Theory of Money},
  editor    = {Friedman, Milton},
  pages     = {25--117},
  publisher = {The University of Chicago Press},
  address   = {Chicago},
  year      = {1956}
}

@article{gomes-mohr-riagosouza-2013,
  author  = {Gomes, Diogo A. and Mohr, Joana and Souza, Rafael Rig{\~a}o},
  title   = {Continuous Time Finite State Mean Field Games},
  journal = {Applied Mathematics \& Optimization},
  volume  = {68},
  number  = {1},
  pages   = {99--143},
  year    = {2013}
}

@article{gomes-mohr-riagosouza,
  author  = {Gomes, Diogo A. and Mohr, Joana and Souza, Rafael Rig{\~a}o},
  title   = {Discrete Time, Finite State Space Mean Field Games},
  journal = {Journal de Math\'ematiques Pures et Appliqu\'ees},
  volume  = {93},
  number  = {3},
  pages   = {308--328},
  year    = {2010}
}
}

\end{document}